\documentclass[11pt]{article}
\usepackage{float}
\usepackage{amsmath}
\usepackage{graphicx}
\usepackage{enumerate}
\usepackage{natbib}
\usepackage{url} 
\usepackage{epstopdf}
\usepackage{epsfig}
\usepackage{multirow}
\usepackage{graphics}
\usepackage{subfigure}
\usepackage{caption}
\usepackage{float}

\newtheorem{thm}{Theorem}

\newtheorem{prop}{Proposition}

\begin{document}

\def\spacingset#1{\renewcommand{\baselinestretch}%
{#1}\small\normalsize} \spacingset{1}


  \title{\bf A Unified Approach for Constructing Confidence Intervals and Hypothesis Tests Using h-function}
  \author{Weizhen Wang
  	 \hspace{.2cm}\\
  		Department of Mathematics and Statistics\\ Wright
  		State University
  		}
\date{March 4, 2021}
  \maketitle


\bigskip
\begin{abstract}
We introduce a general method, named the h-function method, to unify the constructions of level-$\alpha$ exact test and $1-\alpha$ exact confidence interval. Using this method, any confidence interval is improved as follows: i) an approximate interval, including a point estimator, is modified to an exact interval; ii) an exact interval is refined to be an interval that is a subset of the previous one. Two real datasets are used to illustrate the method. 
\end{abstract}

\noindent%

\vfill
\spacingset{1.5} 

\section{Introduction}
It is well known that there is a one-to-one mapping between a family of tests and a confidence set. If one of the two is given in advance, then the other can be derived as follows. 

Let $\Theta$ be the range of a parameter of interest $\theta$ and $S$ be a sample space.
For each $\theta_0\in \Theta$, let $A(\theta_0)$ be the acceptance region of a level-$\alpha$ test of $H_0: \theta= \theta_0$. Then, 
\begin{equation}
\label{ci}
C(\underline{x})=\{\theta_0\in \Theta: \underline{x} \in A(\theta_0)\}.
\end{equation}
is a $1 - \alpha$ confidence set. Conversely, let $C(\underline{x})$ be a $1-\alpha$ confidence set. Define
\begin{equation}
\label{test}
A(\theta_0)=\{\underline{x}\in S: \theta_0\in C(\underline{x})\}.
\end{equation}
Then $A(\theta_0)$ is the acceptance region of a level-$\alpha$ test of $H_0: \theta=\theta_0$.

An interval is typically derived from  the tests since tests are easier to construct. However, solving $C(\underline{x})$ from $A(\theta_0)$ as in (\ref{ci}) is a complicated process. 
The first goal of this paper is to simplify this process.  This is different from the existing methods described in Casella and Berger (2002), where the test and interval are constructed separately. 

The key feature of a confidence interval is that its confidence coefficient, which is defined to be the infimum coverage probability (ICP) over the entire parameter space  (Casella and Berger, 2002), should be no smaller than the nominal level $1-\alpha$. To avoid ambiguity in the future discussion, a 
$1-\alpha$ exact confidence interval means that it has an ICP no less than $1-\alpha$; while a $1-\alpha$ confidence interval only means that the nominal level of the interval is set to be $1-\alpha$ but its ICP can be any number in $[0,1]$. For example, a 95\% Wald internal for a proportion is not necessarily a 95\% exact interval. In fact, a $1-\alpha$ Wald interval is a zero exact interval for any sample size and any $\alpha$ in $[0,1]$ (Brown, Cai and DasGupta, 2001; Agresti, 2013). 

The requirement of an ICP no smaller than $1-\alpha$ is violated for an approximate interval. Thus, the major concern is whether an inferential conclusion drawn from the interval is reliable as its ICP is seldom reported and can be much smaller than the nominal level $1-\alpha$. Huwang (1995) proved this for Wilson interval (1927) even for large samples. On the other hand, an approximate interval is attractive in practice since it is easy to derive.   
Therefore, it is of great interest to build an exact interval through a given approximate interval -- the second goal of the paper.  To the best of our knowledge, limited research has been conducted on this issue. 

One common way to obtain a $1-\alpha$ exact two-sided interval is to take the intersection of two one-sided $1-{\alpha\over 2}$ intervals, for example, Clopper-Pearson interval (1934). Such intervals may be conservative. Hence, it is also of great interest to shrink any given $1-\alpha$ exact interval -- the third goal of the paper. Casella (1986) and Wang (2014) refine exact confidence intervals for a single proportion. However, their methods cannot be applied to the general case when there exist nuisance parameters.

We address the above three problems by introducing an h-function that is closely related to the p-value. 
The computation of p-values was originally done by Arbuthnot (1710) in a study of comparing the  probabilities of male and female births. Laplace addressed the same problem later using binomial distributions (Stigler, 1986). The concept of the p-value was first formally introduced by Pearson (1900) and the use of the p-value in Statistics was popularized by Fisher (1925), where he proposed the level $0.05$ as a limit for statistical significance.
Now, the p-value method is a commonly used method for test construction, where the p-value is treated as a function over the random observation. In this paper, however, the p-value is also treated as a function over the parameter of interest for confidence interval construction. 
It can modify and/or refine any interval
-- a solution to the challenging problem of improving an interval when a nuisance parameter exists.

The paper is organized as follows: In Section 2, we describe the h-function  method that yields both level-$\alpha$ exact tests and $1-\alpha$ exact confidence intervals. Section 3 discusses how to modify any two-sided confidence interval to an exact interval. Section 4 improves any exact two-sided interval. Section 5 modifies any one-sided interval to the smallest one. We provide discussions in Section 6. All proofs are given in the Appendix. 

\section{A general method}

Suppose $\underline{X}$ is observed from a distribution with joint cumulative distribution function (CDF) $F_{(\theta,\underline{\eta})}(\underline{x})$ specified by a parameter vector $(\theta,\underline{\eta})$ in a parameter space $H$. Here $\theta$ is the parameter of interest and $\underline{\eta}$ is the nuisance parameter
vector. The null hypothesis $H_0$ is one of the three forms: $\theta=\theta_0$, $\theta\leq \theta_0$ and $\theta\geq \theta_0$, for a fixed value $\theta_0$. Each form corresponds to two-sided, lower and upper one-sided intervals for $\theta$, respectively. We next introduce the h-function method and illustrate its basic usage through two simple cases.

\subsection{The h-function method}
 A p-value $p(\underline{X})$ is a statistic satisfying $0 \leq p(\underline{x})\leq 1$ for every sample point $\underline{x}$. A small $p(\underline{x})$ provides evidence that $H_0$ is false. Following Casella and Berger (2002),
a p-value is valid if, for every $0 \leq \alpha\leq 1$,
\begin{equation}
\label{pvalue-valid}
\sup_{\theta\in H_0}P_{(\theta,\underline{\eta})} (p(\underline{X})\leq \alpha) \leq \alpha.
\end{equation} 
For simplicity, we drop the subscript $(\theta,\underline{\eta})$ in the future discussion. 
In many cases, a valid p-value $p(\underline{X})$ at $\underline{x}$ can be defined through a test statistic $T(\underline{X})$ using, for example, 
\begin{equation}\label{p-ts-small}
p(\underline{x})=\sup_{\theta\in H_0} P(T(\underline{X})\leq T(\underline{x})),
\end{equation}
when a small value of $T(\underline{X})$ is not in favor of $H_0$. This includes the case that $T(\underline{X})$ is the likelihood ratio test statistic for $H_0$. For any  $p(\underline{X})$ satisfying (\ref{pvalue-valid}), a level-$\alpha$ acceptance region for $H_0$ is equal to  
$\{\underline{x}: p(\underline{x})> \alpha\}.$

Note that $H_0$ is equal to a single point $\theta_0$ or a one-sided interval  with the boundary point $\theta_0$. The p-value $p(\underline{X})$ indeed depends on both $\underline{X}$ and $\theta_0$. 
 This is the key fact for the future theoretical development. Thus, we rewrite the p-value as
\begin{equation}
\label{h-function}
h(\underline{X},\theta_0)\stackrel{def}{=}p(\underline{X})
\end{equation}
and call the left hand side the h-function. We emphasis that the h-function closely relates to but is different from the p-value. The former is a function of both the random observation and  parameter of interest, while the latter is typically treated as a function of the random observation. 
Based on this $h(\underline{X},\theta_0)$, the level-$\alpha$ acceptance region for $H_0$ and $1-\alpha$ exact confidence set for $\theta$ are given by 
\begin{equation}
\label{test-ci}
A(\theta_0)=\{\underline{x}: h(\underline{x},\theta_0)>\alpha\}\,\ \mbox{and}\,\ C(\underline{x})=\{\theta_0: h(\underline{x},\theta_0)>\alpha\},
\end{equation}
respectively.  Both $A(\theta_0)$ and $ C(\underline{x})$ are obtained by solving the same inequality $h(\underline{x},\theta_0)>\alpha$ but in terms of two different arguments $\underline{x}$ and $\theta_0$, respectively. In this sense, the constructions of test and confidence set are unified. We do not obtain $A(\theta_0)$ or $ C(\underline{x})$ from each other following (\ref{ci}) and (\ref{test}). Instead, we use the intermediary h-function in (\ref{h-function}) and name it the h-function method. In fact, Blaker (2000) used  a special $h$ in (\ref{test-ci}) to derive confidence intervals for four discrete distributions of a single parameter. Here, the h-function method is applied to a general case with or without nuisance parameters. In particular, Blaker interval (2000) can be uniformly shortened by this method as shown in Table \ref{table-4}.

It is known that $C(\underline{x})$ may not be an interval. So, we use $\overline{C(\underline{x})}$. Here, 
$\overline{A}$ denotes the smallest simply connected set containing set $A$. Therefore,  $\overline{C(\underline{x})}$ is always an interval and its infimum coverage probability (ICP) over the entire parameter space $H$ is at least $1-\alpha$ because $\overline{C(\underline{x})}$ contains $C(\underline{x})$. i.e.,
$$ICP(\overline{C(\underline{X})})=\inf_{(\theta,\underline{\eta})\in H}P(\theta\in \overline{C(\underline{X})})\geq \inf_{(\theta,\underline{\eta})\in H}P(\theta\in C(\underline{X}))=ICP(C(\underline{X}))
\geq 1-\alpha.$$
The ICP is also called the confidence coefficient (Casella and Berger, 2002).
Throughout the paper, we use $C(\underline{X})$ to denote  $\overline{C(\underline{X})}$ and $\overline{\{\theta_0: h(\underline{x},\theta_0)>\alpha\}}$
if it causes no confusion.

In general, a test statistic may also depend on $\theta_0$, for example, the t-statistic, and thus has a form of $T(\underline{X},\theta_0)$. When a small value of $T(\underline{X},\theta_0)$ is in favor of $H_A$, let $K(\underline{x},\theta_0)=\{\underline{y}:T(\underline{y},\theta_0)\leq T(\underline{x},\theta_0)\}$ be a subset of the sample space, then 
\begin{equation}
\label{h-function-ts-small}
h(\underline{x},\theta_0)
=\sup_{(\theta,\underline{\eta})\in H_0}P(K(\underline{x},\theta_0))
=\left\{ \begin{split}
&\sup_{(\theta,\underline{\eta}) \in H_0}\sum_{\underline{y}\in K(\underline{x},\theta_0)}  f_{(\theta,\underline{\eta})}(\underline{y})\\
&\sup_{(\theta,\underline{\eta}) \in H_0}\int_{K(\underline{x},\theta_0)}  f_{(\theta,\underline{\eta})}(\underline{y})d\underline{y},
\end{split}
\right.
\end{equation}
where $f_{(\theta,\underline{\eta})}$ is either the joint probability mass function (PMF) or probability density function (PDF) of $\underline{X}$.

For illustration of the h-function method, 
we next consider two classic problems: i) estimating the proportion $p$ based on a binomial observation $X\sim Bino(n,p)$;  ii) estimating
the difference of two proportions $d=p_1-p_2$ based on two independent binomials $X\sim Bino(n_1,p_1)$ and $Y \sim Bino(n_2,p_2)$. These are widely used in practice, including clinical trials. The two problems are still open since the best intervals for $p$ and $d$ have not been recognized yet. Let $p_B(x,n,p)$ and $F_B(x,n,p)$ denote the PMF and CDF of $Bino(n,p)$.

\subsection{Two-sided confidence intervals for a proportion}
For a fixed $p_0\in [0,1]$, consider the hypotheses
$H_0: p = p_0\,\ \mbox{vs} \,\ H_A: p\not= p_0.$ Suppose $T_p(x,p_0)$ is a test statistic satisfying 
\begin{equation}\label{p-ci-l-u}
T_p(x,p_0)=T_p(n-x,1-p_0), \,\ \forall x\in [0,n],
\end{equation}
and a small value of $T_p(x,p_0)$ supports $H_A$. Hence, the h-function based on $T_p(x,p_0)$ is
$$h_p(x,p_0)=\sum_{\{y\in [0,n]:T_p(y,p_0)\leq T_p(x,p_0)\}} p_B(y,n,p_0).$$
The acceptance region of level-$\alpha$ test for $H_0$ and $1-\alpha$ exact confidence interval for $p$ are 
$$A_p(p_0)=\{x: h_p(x,p_0)>\alpha\}\,\ \mbox{and}\,\ C_p(x)=\overline{\{p_0: h_p(x,p_0)>\alpha\}}.$$
Now we state a fact that the upper limit $U_p(x)$ of interval $C_p(x)$ can be obtained using the lower limit $L_p(x)$ and vice versa. This simplifies interval construction by half.

\begin{prop}\label{prop-p-l-u} 
	For a test statistic $T_p$ satisfying (\ref{p-ci-l-u}), we have
	\begin{equation}\label{p-l-u}
	U_p(x)=1-L_p(n-x),\,\ \forall x \in [0,n].
	\end{equation}
\end{prop}

Define the first h-function $$h_{p1}(x,p_0)=\min\{2\min\{P_{p_0}(X\leq x), P_{p_0}(X\geq x)\}, 1\}.$$ This yields the two-sided Clopper-Pearson interval (1934) for $p$, denoted by $C_{p1}(X)$.  The interval also satisfies (\ref{p-l-u}) since $h_{p1}(x,p_0)=h_{p1}(n-x,1-p_0)$.

 Define the second h-function based on a test statistic $T_{p2}(x,p_0)$, $$h_{p2}(x,p_0)=\sum_{\{y\in [0,n]:T_{p2}(y,p_0)\leq T_{p2}(x,p_0)\}}p_B(y,n,p_0),$$ where $T_{p2}(x,p_0)=\min\{P_{p_0}(X\leq x), P_{p_0}(X\geq x)\}$.
 This yields  Blaker interval (2000), denoted by $C_{p2}(X)$, for $p$. Since $T_{p2}(x,p_0)$ satisfies (\ref{p-ci-l-u}), interval $C_{p2}(X)$ satisfies (\ref{p-l-u}). As discussed in Casella and Berger (2002) and Agresti (2013), $C_{p2}(X)$ has a nesting property: an interval $C_{p2}(X)$ with a higher confidence level always contains the one with a lower level. Also, $C_{p2}(X)$ is a subset of $C_{p1}(X)$. But $C_{p2}(X)$ can still be shortened uniformly as shown in Section 4.

The third h-function based on the likelihood ratio test statistic $T_{p3}(x,p_0)$ is 
$$h_{p3}(x,p_0)=\sum_{\{y\in [0,n]: T_{p3}(y,p_0)\leq T_{p3}(x,p_0)\}} p_B(y,n,p_0),$$
where $T_{p3}(x,p_0)=({p_0\over \hat{p}})^x({1-p_0\over 1-\hat{p}})^{n-x}$ for $\hat{p}=x/n$. This generates a confidence interval $C_{p3}(x)$.
The upper limit of $C_{p3}(x)$ can also be determined by its lower limit following (\ref{p-l-u}) because $T_{p3}(x,p_0)$ satisfies (\ref{p-ci-l-u}).

We report in Table~\ref{table-1} the infimum coverage probability (ICP) and total interval length (TIL) over all $n+1$ sample points for the above intervals $C_{pi}$ for $i=1,2,3$ when $n$ is from 16 to 100. Their ICPs are all no smaller than 0.95. Due to Wang (2007), the ICP of an interval $C(X)$ for $p$ with a nondecreasing lower confidence limit $L(X)$ is achieved at one of the values $(L(x))^-$ for $x=1,...,n$, where $a^-$ denotes the left limit of $y$ when $y$ approaches $a$. Thus, the ICP of $C(X)$ can be computed precisely.
Interval $C_{p2}$ is the winner among the three due to its small TIL. In Section 4, we see that these three intervals are uniformly shorter, as shown in Table~\ref{table-4}, than their modified intervals $C^M_{pi}$. For comparison purpose, the interval  $C_{p4}$ in Wang (2014) is also reported. This interval, derived by an iterative algorithm, is admissible and is shortest among the four intervals. Here, a $1-\alpha$ exact confidence interval $C(\underline{X})$ is admissible if any interval $C'(\underline{X})$, which is contained in but not equal to $C(\underline{X})$, has an ICP less than $1-\alpha$.

\begin{table}[h]
	\begin{center}
		\renewcommand{\arraystretch}{0.9}
		\caption{The infimum coverage probability (ICP) and total interval length (TIL) of intervals over all sample points for eight exact confidence intervals for $p$: $C_{pi}$ and $C^M_{pi}$ for $i=1$ (Clopper-Pearson), 2 (Blaker), 3 (Likelihood-ratio-test interval), 4 (Wang),  when $1-\alpha=0.95$ and $n$ varies. The smallest TIL for each $n$ is marked by `*'.} 
		\begin{tabular}{l|lll|lll}
			$n$	  &          & ICP      & TIL       &           & ICP      & TIL \\	\hline	 
			$16$  &	$C_{p1}$ & 0.9578& 6.9380& $C^M_{p1}$& 0.9500& 6.4978*\\
			& $C_{p2}$ & 0.9500& 6.5043& $C^M_{p2}$& 0.9500& 6.4978*\\
			& $C_{p3}$ & 0.9500& 6.6115& $C^M_{p3}$& 0.9500& 6.5342\\
			& $C_{p4}$ & 0.9500& 6.4978*&$C^M_{p4}$& 0.9500& 6.4978*\\ \hline
			$30$  & $C_{p1}$ & 0.9505& 9.2705& $C^M_{p1}$& 0.9500& 8.7784\\
			& $C_{p2}$ & 0.9500& 8.7814& $C^M_{p2}$& 0.9500& 8.7770\\
			& $C_{p3}$ & 0.9500& 8.8742& $C^M_{p3}$& 0.9500& 8.8089\\
			& $C_{p4}$ & 0.9500& 8.7726*&$C^M_{p4}$& 0.9500& 8.7726*\\ \hline
			$100$  & $C_{p1}$ & 0.9503& 16.3057&$C^M_{p1}$& 0.9500& 15.8214\\
			& $C_{p2}$ & 0.9500& 15.8243&$C^M_{p2}$& 0.9500& 15.8176\\
			& $C_{p3}$ & 0.9500& 15.8803&$C^M_{p3}$& 0.9500& 15.8402\\
			& $C_{p4}$ & 0.9500& 15.8146*& $C^M_{p4}$& 0.9500& 15.8146*
		\end{tabular}
		\label{table-1}
	\end{center}
\end{table}

\subsection{Confidence intervals for the difference of two proportions}

 There are three commonly used measurements, the difference $d=p_1-p_2$, the relative risk and the odds ratio, for comparison of the two proportions. Here we focus on the difference.
Consider the following hypotheses for any fixed $d_0\in [-1,1]$:
$H_0: d = d_0\,\ \mbox{vs} \,\ H_A: d\not= d_0.$ Under $H_0$, $p_1=d_0+p_2$ for $p_2\in D(d_0)$, 
where
\begin{eqnarray*}
	D(d_0)=\left\{
	\begin{split}
		& [0,1-d_0], \,\ \mbox{if} \,\ d_0\in [0,1];\\
		&  [-d_0,1], \,\ \mbox{if} \,\ d_0\in [-1,0).
	\end{split}
	\right.
\end{eqnarray*}

Suppose $T_d(x,y,d_0)$ is a test statistic satisfying 
\begin{equation}\label{diff-td}
T_d(x,y,d_0)=T_d(n_1-x,n_2-y,-d_0), \,\ \forall (x,y)\in S_d,
\end{equation}
for $S_d=[0,n_1]\times[0,n_2]$ and a small value of $T_d(x,y,d_0)$ supports $H_A$. The h-function based on $T_d(x,y,d_0)$ is
\begin{equation}
\label{h-diff}
h_d(x,y,d_0)=\sup_{p_2\in D(d_0)} \sum_{ \{(u,v)\in S_d: T_d(u,v,d_0)\leq T_d(x,y,d_0)\}}p_B(u,n_1,p_2+d_0) p_B(v,n_2,p_2).
\end{equation}
The acceptance region of level-$\alpha$ test for $H_0$ and $1-\alpha$ exact confidence interval for $d$ are 
\begin{equation}
\label{test-ci-diff}
A_d(d_0)=\{(x,y): h_d(x,y,d_0)>\alpha\}\,\ \mbox{and}\,\ C_d(x,y)=\overline{\{d_0: h_d(x,y,d_0)>\alpha\}}.
\end{equation}
Similar to Proposition~\ref{prop-p-l-u}, we state a fact that determines the upper limit $U_d(x,y)$ of $C_d(x,y)$ from its lower limit $L_d(x,y)$ and vice versa. This again simplifies  interval calculation.

\begin{prop}\label{prop-diff-l-u} 
	For a test statistic $T_d$ satisfying (\ref{diff-td}), we have
	\begin{equation}\label{diff-ci-l-u}
	U_d(x,y)=-L_d(n_1-x,n_2-y),\,\ \forall (x,y) \in S_d.
	\end{equation}
\end{prop}

First, introduce the likelihood ratio test statistic:
\begin{eqnarray*}
	T_{d1}(x,y,d_0)
	=({\hat{p}_{1d}\over \hat{p}_1})^x ({1-\hat{p}_{1d} \over 1-\hat{p}_1})^{n_1-x}  ({\hat{p}_{2d}\over \hat{p}_2})^y ({1-\hat{p}_{2d} \over 1-\hat{p}_{2}})^{n_2-y},
\end{eqnarray*}
where $$\hat{p}_{2d}(x,y,d_0)=arg\max_{p_2\in D(d_0)}p_B(x, n_1,p_2+d_0)p_B(y,n_2,p_2), \,\ $$ $$ \hat{p}_{1d}(x,y,d_0)=\hat{p}_{2d}(x,y,d_0)+d_0, \,\ \hat{p}_1=x/n_1, \,\ \hat{p}_2=y/n_2.$$
It can be shown that $T_{d1}$ satisfies (\ref{diff-td}). The h-function $h_{d1}(x,y,d_0)$ based on $T_{d1}$ 
 follows (\ref{h-diff}).
Then, the acceptance region of level-$\alpha$ test and $1-\alpha$ confidence interval for $d$ are given in (\ref{test-ci-diff}) and are denoted by $A_{d1}(d_0)$ and $C_{d1}(x,y)$, respectively. The upper limit of $C_{d1}(x,y)$ is determined by the lower limit through (\ref{diff-ci-l-u}).

Secondly, the score test statistic, see Agresti and Min (2001) and Fay (2010), is  
$$T_{d2}(x,y,d_0)={-|\hat{p}_1-\hat{p}_2-d_0| \over\sqrt{{\hat{p}_{1d}(x,y,d_0)(1-\hat{p}_{1d}(x,y,d_0))\over n_1} +{\hat{p}_{2d}(d_0)(1-\hat{p}_{2d}(x,y,d_0))\over n_2} } }.$$
When $(x,y,d_0)=(n_1,0,1)$ or $(0,n_2,-1)$, the above ratio is 0/0. So it is defined to be 0.
We obtain the h-function $h_{d2}(x,y,d_0)$ following (\ref{h-diff}) and then derive the acceptance region of level-$\alpha$ test $A_{d2}(d_0)$ and $1-\alpha$ exact confidence interval  $C_{d2}(x,y)$ following  (\ref{test-ci-diff}). Since $T_{d2}$ satisfies (\ref{diff-td}), the  
upper limit of $C_{d2}(x,y)$ is determined by the lower limit through (\ref{diff-ci-l-u}).

Table~\ref{table-2} reports the ICP and TIL over
all $(n_1+1)(n_2+1)$ sample points in $S_d$ for $C_{d1}$ and $C_{d2}$ and other intervals to be discussed later when $1-\alpha=0.95$ and $(n_1,n_2)$ varies. The ICPs of  $C_{d1}$ and $C_{d2}$ are at least 0.95. Each value is computed by a grid search: select $201^2$ pairs of $(p_1,p_2)$, where both $p_1$ and $p_2$ are the multiples of 0.005; compute the coverage probability at each pair; use the minimum of these coverage probabilities as the ICP. 
\begin{table}[h]
	\begin{center}
		\renewcommand{\arraystretch}{0.9}
		\caption{The ICPs and TILs over all sample points in $S_d$ for fifteen intervals $(C_{di}$, $C^M_{di}$, $C^{M\infty}_{di})$ for $i=1$ (LRT, exact), $2$ (Score, exact), $3$ (Wald, approximate), $4$ (Maximum likelihood estimator (MLE) as an interval) and $5$ (Wang, exact) when $1-\alpha=0.95$ and $(n_1,n_2)$ varies. The smallest TIL for exact intervals is marked by `*' for each $(n_1,n_2)$.} 
		\resizebox{1\textwidth}{0.2\textheight}{
			\begin{tabular}{l|ccc|ccc|ccc}
				$(n_1,n_2)$&			& ICP & TIL & & ICP & TIL &  & ICP & TIL  \\ \hline
				(8,10)	   & $C_{d1}$ &  0.9500 & 77.2224 & $C^M_{d1}$ & 0.9500 & 75.7339 & $C^{M\infty}_{d1}$ & 0.9500  & 74.9249\\
				& $C_{d2}$ &  0.9500 & 73.3681 & $C^M_{d2}$ & 0.9500 & 72.9133 & $C^{M\infty}_{d2}$ & 0.9500  & 72.6113\\
				& $C_{d3}$ &  0 & 67.8756 & $C^M_{d3}$ & 0.9503 & 96.3142 & $C^{M\infty}_{d3}$ & 0.9500  & 83.2143\\
				& $C_{d4}$ &  0 & 0 & $C^M_{d4}$ & 0.9500 & 76.8064 & $C^{M\infty}_{d4}$ & 0.9500  & 76.2304\\
				& $C_{d5}$ &  0.9515 & 76.9506 & $C^M_{d5}$ & 0.9500 & 73.0494 & $C^{M\infty}_{d5}$ & 0.9500  & 72.4728*\\
				\hline
				(10,15)	   & $C_{d1}$ &  0.9500 & 119.9192 & $C^M_{d1}$ & 0.9500 & 117.2285 & $C^{M\infty}_{d1}$ & 0.9500  & 115.6740\\
				& $C_{d2}$ &  0.9500 & 113.3737 & $C^M_{d2}$ & 0.9500 & 112.1987 & $C^{M\infty}_{d2}$ & 0.9500  & 111.5613*\\
				& $C_{d3}$ &  0 & 106.2471 & $C^M_{d3}$ & 0.9500 & 148.7108 & $C^{M\infty}_{d3}$ & 0.9500  & 132.0608\\
				& $C_{d4}$ &  0 & 0 & $C^M_{d4}$ & 0.9500 & 120.7789 & $C^{M\infty}_{d4}$ & 0.9500  & 120.7007\\
				& $C_{d5}$ &  0.9515 & 116.8048 & $C^M_{d5}$ & 0.9500 & 112.6569 & $C^{M\infty}_{d5}$ & 0.9500  & 111.7894\\
				\hline
				(23,32)	   & $C_{d1}$ &  0.9500 & 361.8775 & $C^M_{d1}$ & 0.9500 & 358.1373 & $C^{M\infty}_{d1}$ & 0.9500  & 355.8710\\
				& $C_{d2}$ &  0.9500 & 346.4825 & $C^M_{d2}$ & 0.9500 & 344.3728 & $C^{M\infty}_{d2}$ & 0.9500 & 342.6230\\
				& $C_{d3}$ &  0 & 332.3962 & $C^M_{d3}$ & 0.9500 & 399.0738 & $C^{M\infty}_{d3}$ & 0.9500  & 380.5928\\
				& $C_{d4}$ &  0 & 0 & $C^M_{d4}$ & 0.9500 & 372.2596 & $C^{M\infty}_{d4}$ & 0.9500  & 370.0785\\
				& $C_{d5}$ &  0.9503 & 347.4601 & $C^M_{d5}$ & 0.9500 & 342.6516 & $C^{M\infty}_{d5}$ & 0.9500  & 341.4697*\\
			\end{tabular}
			\label{table-2}
		}
	\end{center}
\end{table}

\section{Modifying any given two-sided confidence interval}

For convenience, we consider the closed interval  $C_0(\underline{X})=[L_0(\underline{X}),U_0(\underline{X})]$ for $\theta$ if the confidence limits are not infinity. 
For a given value $\theta_0$, consider 
\begin{equation}\label{hy-ci-test}
H_0: \theta=\theta_0\,\ \mbox{vs}\,\ H_A: \theta\not=\theta_0.
\end{equation}
Define a test statistic through interval $C_0(\underline{X})$ 
\begin{equation}\label{ts-t2}
T_2(\underline{X},\theta_0)=\min\{\theta_0-L_0(\underline{X}), U_0(\underline{X})-\theta_0 \}.
\end{equation}  
Clearly, a small  $T_2(\underline{x},\theta_0)$ provides strong evidence of establishing $H_A$,
and 
\begin{equation}
\label{t2-c0}
T_2(\underline{x},\theta_0)\geq 0 \,\ \mbox{if and only if} \,\ \theta_0\in C_0(\underline{x}).
\end{equation}
The h-function based on $T_2(\underline{x},\theta_0)$ is
\begin{equation}
\label{h-2}
h_2(\underline{x},\theta_0)=\sup_{H_0}P(T_2(\underline{X},\theta_0)\leq T_2(\underline{x},\theta_0)).
\end{equation}
Following (\ref{test-ci}), the level-$\alpha$ acceptance region for $H_0$ and $1-\alpha$ exact two-sided confidence interval for $\theta$ are 
\begin{equation}
\label{test-ci-2}
A_2(\theta_0)=\{\underline{x}: h_2(\underline{x},\theta_0)>\alpha\}\,\ \mbox{and}\,\ C_0^M(\underline{x})=\overline{\{\theta_0: h_2(\underline{x},\theta_0)>\alpha\}}.
\end{equation}
In the rest of the paper, $\overline{A}$ denotes the smallest closed simply connected set that contains set $A$.
The superscript ``M'' refers a modification. In general, $C_0^M(\underline{x})\not = C_0(\underline{X})$.
For a positive integer $k$, $C_0^{Mk}(\underline{x})$ is the resultant interval after the modification process of (\ref{ts-t2}) through (\ref{test-ci-2}) is applied to $C_0(\underline{X})$ for $k$ consecutive times.
For example, when  $k=1$, $C_0^{M1}(\underline{x})=C_0^M(\underline{x})$, and when $k=2$, $C_0^{M2}(\underline{x})=(C_0^M)^M(\underline{x}).$ 

\begin{thm}
	\label{thm-ci-modi}
	For any given confidence interval $C_0(\underline{X})$, consider the hypotheses (\ref{hy-ci-test}).
	
	i) The h-function $h_2(\underline{x},\theta_0)$ in (\ref{h-2}) is a valid p-value for statistic $T_2(\underline{X},\theta_0)$ in (\ref{ts-t2}) at observation $\underline{x}$.  
	
	ii) Interval
	$C_0^M(\underline{X})$ given in (\ref{test-ci-2}) is of level $1-\alpha$. i.e., $ICP(C_0^M)\geq 1-\alpha$. 
\end{thm}

The theorem modifies a confidence interval $C_0(\underline{X})$ of any level,  including a point estimator (a zero exact confidence interval), to a $1-\alpha$ exact confidence interval $C_0^M(\underline{X})$. Also, when using $C_0^M(\underline{X})$ to conduct a level-$\alpha$ test for $H_0$ one simply rejects $H_0$ if $\theta_0$ falls outside the interval but is unable to report a p-value. Now the p-value is given by $h_2(\underline{x},\theta_0)$ in (\ref{h-2}). The test of rejecting $H_0$ if this p-value is no larger than $\alpha$ is an exact test of level $\alpha$; On the contrary, the test of rejecting $H_0$ if $\theta_0$ falls outside $C_0(\underline{X})$ is not of level $\alpha$ if  $C_0(\underline{X})$ is not a $1-\alpha$ exact interval.
\medskip

{\bf Example 1}. (more confidence intervals for a proportion) Interval estimation of $p$ based on $X\sim Bino(n,p)$ is one of the basic inferential problems.  Asymptotic intervals play an important role but may not be reliable. We apply the modification process of (\ref{ts-t2}) through (\ref{test-ci-2}) to several commonly used intervals to obtain exact intervals.

Let $z_{\alpha}$ denote the upper $\alpha$th percentile of the standard normal distribution. It is known that Wald interval
$$C_{p5}(X)=[ \hat{p} \mp z_{{\alpha \over 2}}\sqrt{{\hat{p}(1-\hat{p})\over n}}]$$
has an ICP zero for any $n$ and $\alpha\in (0,1)$ because it shrinks to a point when $X=0$ or $n$ (Agresti, 2013). We now modify it to a $1-\alpha$ exact interval using the modification process. Table~\ref{table-3} contains 95\% $C_{p5}(X)$,
$C_{p5}^M(X)$, $C_{p5}^{M22}(X)$ and others
when $n=16$. For example,
$C_{p5}(X)$ has $ICP=0$, $TIL=6.0559$ and $C_{p5}(3)=[-0.0038,0.3788]$. In order to assure the computed ICP for an exact interval to be at least $1-\alpha$, the lower limits for a confidence interval round down but the upper limits round up at the fourth decimal place.
Although  $C_{p5}(X)$ has an zero ICP, its modification $C_{p5}^M(X)$ has a correct ICP 95\%. In return,  $C_{p5}^M(X)$ becomes wider and does not  make the lower limits at $x=1,2,3$ larger than zero. This is due to the negative values of the lower limits of $C_{p5}(X)$ at $x=1,2,3$. We further apply the modification process to $C_{p5}^M(X)$ for more times. The 22nd modification 
$C_{p5}^{M22}(X)$ has an ICP 95\% and is a subset of 
$C_{p5}^M(X)$. It cannot be shortened by the modification process any further, however, can be shortened at $x=1,2,3$ following Casella (1986) or Wang (2014).  

We also report, in Table~\ref{table-3}, the result of modifying the 95\% Wilson interval $C_{p6}(X)$ (1927). Huwang (1995) showed that its ICP is much lower than the nominal level $1-\alpha$ for any $n$. However, the modified $C_{p6}^M(X)$ has a correct ICP 95\% and is admissible here. To see the admissibility, we apply the refining approach in Wang (2014) to $C_{p6}^M(X)$ and find that the refinement of $C_{p6}^M(X)$ is equal to  $C_{p6}^M(X)$. Since Wang's approach produces admissible intervals, we conclude $C_{p6}^M(X)$ admissible when $(n,1-\alpha)=(16,0.95)$. 

\begin{table}[H]
	\begin{center}
		\renewcommand{\arraystretch}{1}
		\caption{ 95\% Wald interval $C_{p5}$, the modification
			$C_{p5}^M$, the 22nd modification $C_{p5}^{M22}$; 95\% Wilson interval  $C_{p6}$, 
			$C_{p6}^M$; the sample-proportion estimator $C_{p7}$,  
			$C_{p7}^M$; another point estimator $C_{p8}$, 
			$C_{p8}^M$, $C_{p8}^{M17}$; their ICPs and TILs 
			when $n=16$. For each interval, the first line is the lower limits and the second line is the upper limits. } 
		\resizebox{1\textwidth}{0.42\textheight}{
			\begin{tabular}{l|rrrrrrrrrr}
				$X$ &   0   & 1     & 2     & 3     & 4     & 5     & 6     & 7   & 8 \\	\hline
				$C_{p5}$	& 0.0000&-0.0562&-0.0371&-0.0038& 0.0378& 0.0853& 0.1377& 0.1944& 0.2550\\ 
				ICP=	0\%	& 0.0000& 0.1812& 0.2871& 0.3788& 0.4622& 0.5397& 0.6123& 0.6806& 0.7450\\ 
				$C_{p5}^M$	& 0.0000& 0.0000& 0.0000& 0.0000& 0.0189& 0.0426& 0.0688& 0.0972& 0.1275\\
				ICP=	95\%& 0.1709& 0.2674& 0.3522& 0.4295& 0.5000& 0.5805& 0.6626& 0.8403& 0.8725\\ 
				$C_{p5}^{M22}(=C_{p5}^{M\infty})$  & 0.0000& 0.0000& 0.0000& 0.0000& 0.0902& 0.1321& 0.1708& 0.1708& 0.1708\\
				ICP=    95\%& 0.1709& 0.2674& 0.3522& 0.4295& 0.5000& 0.5706& 0.6479& 0.8292& 0.8292\\ \hline
				$C_{p6}$    & 0.0000& 0.0111& 0.0349& 0.0659& 0.1018& 0.1416& 0.1848& 0.2309& 0.2799\\
				ICP=83.62\% & 0.1937& 0.2833& 0.3603& 0.4301& 0.4950& 0.5560& 0.6136& 0.6683& 0.7201\\
				$C_{p6}^M(=C_{p6}^{M\infty})$& 0.0000& 0.0032& 0.0226& 0.0531& 0.0902& 0.1321& 0.1777& 0.2122& 0.2719\\
				ICP=95\%& 0.2123& 0.3076& 0.3734& 0.4371& 0.5000& 0.5630& 0.6267& 0.6925& 0.7281\\ \hline
				$C_{p7}$    & 0.0000& 0.0625& 0.1250& 0.1875& 0.2500& 0.3125& 0.3750& 0.4375& 0.5000\\
				ICP=0\% & 0.0000& 0.0625& 0.1250& 0.1875& 0.2500& 0.3125& 0.3750& 0.4375& 0.5000\\
				$C_{p7}^M(=C_{p7}^{M\infty})$& 0.0000& 0.0032& 0.0226& 0.0531& 0.0902& 0.1321& 0.1777& 0.2187& 0.2719\\
				ICP=95\%& 0.2188& 0.3125& 0.3750& 0.4375& 0.5000& 0.5625& 0.6250& 0.6875& 0.7281\\ \hline
				$C_{p8}$  & 0.0000& 0.0500& 0.1250& 0.2000& 0.2500& 0.3125& 0.3750& 0.4375& 0.5000\\    
				ICP=0\% & 0.0000& 0.0500& 0.1250& 0.2000& 0.2500& 0.3125& 0.3750& 0.4375& 0.5000\\  
				$C_{p8}^{M}$  & 0.0000& 0.0032& 0.0226& 0.0531& 0.0902& 0.1321& 0.1777& 0.2187& 0.2719\\ 
				ICP=95\%  & 0.2188& 0.3063& 0.3750& 0.4438& 0.5000& 0.5585& 0.6250& 0.6938& 0.7281\\
				$C_{p8}^{M17}(=C_{p8}^{M\infty})$& 0.0000& 0.0032& 0.0226& 0.0531& 0.0902& 0.1321& 0.1777& 0.2187& 0.2719\\
				ICP=95\%  & 0.2188& 0.3063& 0.3750& 0.4416& 0.5000& 0.5585& 0.6250& 0.6938& 0.7281\\
				\hline
				$X$ &   9   & 10    & 11    & 12    & 13    & 14    & 15    & 16    & \\
				\hline
				$C_{p5}$    & 0.3194& 0.3877& 0.4603& 0.5378& 0.6212& 0.7129& 0.8188& 1.0000&\\
				TIL=6.0559    & 0.8056& 0.8623& 0.9147& 0.9622& 1.0038& 1.0371& 1.0562& 1.0000&\\
				$C_{p5}^M$  & 0.1597& 0.3374& 0.4195& 0.5000& 0.5705& 0.6478& 0.7326& 0.8291&\\
				TIL=7.8957    & 0.9028& 0.9312& 0.9574& 0.9811& 1.0000& 1.0000& 1.0000& 1.0000&\\
				$C_{p5}^{M22}$  & 0.1708& 0.3521& 0.4294& 0.5000& 0.5705& 0.6478& 0.7326& 0.8291&\\
				TIL=7.0650    & 0.8292& 0.8292& 0.8679& 0.9098& 1.0000& 1.0000& 1.0000& 1.0000&\\ \hline
				$C_{p6}$    & 0.3317& 0.3864& 0.4440& 0.5050& 0.5699& 0.6397& 0.7167& 0.8063&\\
				TIL=6.0974    & 0.7691& 0.8152& 0.8584& 0.8982& 0.9341& 0.9651& 0.9889& 1.0000&\\
				$C^M_{p6}$    & 0.3075& 0.3733& 0.4370& 0.5000& 0.5629& 0.6266& 0.6924& 0.7877& \\
				TIL=6.4978    & 0.7878& 0.8223& 0.8679& 0.9098& 0.9469& 0.9774& 0.9968& 1.0000&\\ \hline
				$C_{p7}$    & 0.5625& 0.6250& 0.6875& 0.7500& 0.8125& 0.8750& 0.9375& 1.0000& \\
				TIL=0              & 0.5625& 0.6250& 0.6875& 0.7500& 0.8125& 0.8750& 0.9375& 1.0000&\\ 
				$C^M_{p7}$    & 0.3125& 0.3750& 0.4375& 0.5000& 0.5625& 0.6250& 0.6875& 0.7812& \\ 
				TIL=6.4978    & 0.7813& 0.8223& 0.8679& 0.9098& 0.9469& 0.9774& 0.9968& 1.0000&\\ \hline
				$C_{p8}$      & 0.5625& 0.6250& 0.6875& 0.7500& 0.8000& 0.8750& 0.9500& 1.0000& \\
				TIL=0         & 0.5625& 0.6250& 0.6875& 0.7500& 0.8000& 0.8750& 0.9500& 1.0000&\\
				$C^{M}_{p8}$ & 0.3062& 0.3750& 0.4415& 0.5000& 0.5562& 0.6250& 0.6937& 0.7812&\\
				TIL=6.5021    & 0.7813& 0.8223& 0.8679& 0.9098& 0.9469& 0.9774& 0.9968& 1.0000&\\
				$C_{p8}^{M17}$& 0.3062& 0.3750& 0.4415& 0.5000& 0.5584& 0.6250& 0.6937& 0.7812&\\
				TIL=6.4978      & 0.7813& 0.8223& 0.8679& 0.9098& 0.9469& 0.9774& 0.9968& 1.0000&\\ \hline         
		\end{tabular}}
		\label{table-3}
	\end{center}
\end{table}

Here is an interesting part. Can we apply the modification process to a point estimator to obtain a $1-\alpha$ exact interval? The answer is yes. Table~\ref{table-3} reports the sample proportion $C_{p7}(X)={X\over n}$ (a zero-length interval estimator of level zero)  and its modification $C_{p7}^M(X)$. The latter turns out to be admissible following the same argument in the previous paragraph. i.e., the refinement of $C_{p7}^M(X)$ from Wang's approach is equal to $C_{p7}^M(X)$. We try another (arbitrary) point estimator $C_{p8}(X)$, which satisfies (\ref{p-l-u}) and has a nondecreasing confidence limits, see Table~\ref{table-3}. The modified $C_{p8}^M(X)$ has an ICP 95\% and can be shortened further because $C_{p8}^{M17}(X)$ is a subset of $C_{p8}^M(X)$. 
~\raisebox{.5ex}{\fbox{}} 

\medskip

{\bf Example 2}. (more confidence intervals for the difference of two proportions) We now modify two approximate intervals for $d=p_1-p_2$ based on two independent binomials $X$ and $Y$ in Section 2.3. The first is Wald interval for $d$:
$$C_{d3}(X,Y)=[\hat{p}_1-\hat{p}_2\mp z_{\alpha\over 2}\sqrt{{\hat{p}_1(1-\hat{p}_1)\over n_1}+{\hat{p}_2(1-\hat{p}_2)\over n_2}}].$$
This interval has a zero ICP for any $n_1$ and $n_2$ and $\alpha$ (Wang and Zhang, 2014). The second is the MLE of $d$ that is treated as a confidence interval $$C_{d4}(X,Y)=[\hat{p}_1-\hat{p}_2\mp 0].$$ Clearly, this interval has a zero ICP and a zero length. 
Theorem~\ref{thm-ci-modi} is applied to these two approximate intervals to obtain two $1-\alpha$ exact confidence intervals $C^M_{d3}(X,Y)$ and $C^M_{d4}(X,Y)$. The modified intervals have a larger TIL over all sample points than the original intervals because they have a correct ICP that is at least $1-\alpha$. Table~\ref{table-2} contains several examples. Theorem~\ref{thm-ci-modi} is also applied to three exact intervals, $C_{d1}(X,Y)$ and $C_{d2}(X,Y)$, discussed in Section 2.3 and the 
two-one-sided interval in Wang (2010), denoted by $C_{d5}(X,Y)$. This time all TIL's of modified intervals decrease rather than increase. See Table~\ref{table-2}.
~\raisebox{.5ex}{\fbox{}} 
\medskip

{\bf Example 3} (the z-interval and its variants) In this example, we investigate the effect of choosing a point estimator on the modified confidence interval.
Suppose $\{X_1,...,X_n\}$ is a random sample from a normal population $N(\mu,\sigma^2)$ for known $\sigma^2$. Confidence intervals for $\mu$ based on the minimum sufficient statistic $\bar{X}$, including the z-interval, are of interest. Consider $H_0: \mu=\mu_0$ vs $H_A: \mu\not= \mu_0$. If $\mu$ is estimated by a point estimator $a\bar{X}+b$ for two given constants $a>0$ and $b$, then, following (\ref{ts-t2}),
define a test statistic
$$T_{zab}(\bar{X},\mu_0)=\min\{\mu_0-a\bar{X}-b, a\bar{X}+b-\mu_0\} =-|\mu_0-a\bar{X}-b|$$
for $H_0$. This includes $T_{z10}$ as a special case for $(a,b)=(1,0)$, where $\mu$ is estimated by the uniformly minimum variance unbiased estimator (UMVUE) $\bar{X}$.
The h-function is
$$h_{zab}(\bar{x},\mu_0)=P(T_{zab}(\bar{X},\mu_0)\leq T_{zab}(\bar{x},\mu_0)).$$
Note $\bar{X}=\mu_0+{\sigma\over \sqrt{n}}Z$ for $Z\sim N(0,1)$ with a CDF $F_N(z)$. 
Then,
\begin{eqnarray*}
	&&h_{zab}(\bar{x},\mu_0)= P(-|\mu_0-a(\mu_0+{\sigma\over \sqrt{n}}Z)-b| \leq  T_{zab}(\bar{x},\mu_0))\\
	&=& \left\{
	\begin{split}
		&1-F_N({\sqrt{n}\over \sigma} [{(2-a)\mu_0-2b-a\bar{x}\over a}])+F_N({\sqrt{n}\over \sigma} (\bar{x}-\mu_0)), \,\ \mbox{if}\,\ \mu_0\geq a\bar{x}+b\\	
		&1-F_N({\sqrt{n}\over \sigma} (\bar{x}-\mu_0))+F_N({\sqrt{n}\over \sigma} [{(2-a)\mu_0-2b-a\bar{x}\over a}]), \,\ \mbox{if}\,\ \mu_0\leq a\bar{x}+b.\\
	\end{split}\right.	
\end{eqnarray*}
This function achieves its maximum 1 at $\mu_0=a\bar{x}+b$.
The $1-\alpha$ exact interval for $\mu$ is 
$$C_{zab}(\bar{x})=\overline{\{\mu_0: h_{zab}(\bar{x},\mu_0)>\alpha\}}\stackrel{denoted \,\ by}{=}[L_{zab}(\bar{x}), U_{zab}(\bar{x})].$$

i) If $a>2$, $C_{zab}(\bar{x})=(-\infty,+\infty)$ since $h_{zab}(\bar{x},\mu_0)$ approaches 1 as $\mu_0$ goes to $\pm \infty$. 

ii) If $a<2$, $C_{zab}(\bar{x})$ is a finite interval  since $h_{zab}(\bar{x},\mu_0)$, as a function of $\mu_0$ for a fixed $\bar{x}$, is unimodal with a mode at $\mu_0=a\bar{x}+b$ and approaches 0 as $\mu_0$ goes to $\pm \infty$.  The confidence limits $L_{zab}(\bar{x})$ and $U_{zab}(\bar{x})$ are the solutions of the equation $h_{zab}(\bar{x},\mu_0)=\alpha$ on intervals $(-\infty, a\bar{x}+b]$ and $[a\bar{x}+b,+\infty)$, respectively. The two solutions  can be found numerically. For the case of $(a,b)=(1,0)$, $h_{z10}(\bar{x},\mu_0)=1-F_N({\sqrt{n}\over \sigma}|\mu_0-\bar{x}|)+F_N(-{\sqrt{n}\over \sigma}|\mu_0-\bar{x}|).$ Then, 
$L_{z10}(\bar{x})=\bar{x}- z_{{\alpha\over 2}}{\sigma \over \sqrt{n}}$ and $U_{z10}(\bar{x})=\bar{x}+ z_{{\alpha\over 2}}{\sigma \over \sqrt{n}}$. i.e., $C_{z10}(\bar{X})$ is 
the z-interval, which is the shortest (Casella and Berger, 2002). 

iii) If $a=2$, then
\begin{eqnarray*}
	h_{z2b}(\bar{x},\mu_0)
	= \left\{
	\begin{split}
		&1-F_N(-{\sqrt{n}\over \sigma} (b+\bar{x}))+F_N({\sqrt{n}\over \sigma} (\bar{x}-\mu_0)), \,\ \mbox{if}\,\ \mu_0\geq 2\bar{x}+b\\	
		&1-F_N({\sqrt{n}\over \sigma} (\bar{x}-\mu_0))+F_N(-{\sqrt{n}\over \sigma} (b+\bar{x})), \,\ \mbox{if}\,\ \mu_0\leq 2\bar{x}+b.\\
	\end{split}\right.	
\end{eqnarray*}
This, as a function of $\mu_0$, is unimodal with a mode at $\mu_0=2\bar{x}+b$. Thus, $U_{z2b}(\bar{x})$ and $L_{z2b}(\bar{x})$ are the solutions of $h_{z2b}(\bar{x},\mu_0)=\alpha$ on intervals $[2\bar{x}+b,+\infty)$ and $(-\infty, 2\bar{x}+b]$, respectively, if a solution exists. Algebraic calculations yield the following three results: 

iii-1) If $\bar{x}> z_{\alpha}{\sigma\over \sqrt{n}}-b$, then  $$L_{z2b}(\bar{x})= \bar{x}-{ \sigma\over \sqrt{n}} F^{-1}_N(1-\alpha+F_N(-{\sqrt{n}\over \sigma}(b+\bar{x})))
\,\ \mbox{and}\,\ U_{z2b}(\bar{x})=+\infty;$$ 

iii-2) If $ -z_{\alpha}{\sigma\over \sqrt{n}}-b\leq \bar{x}\leq  z_{\alpha}{\sigma\over \sqrt{n}}-b$, then  $$L_{z2b}(\bar{x})= -\infty\,\ \mbox{and}\,\ U_{z2b}(\bar{x}) =+\infty;$$

iii-3) If $ \bar{x}<-z_{\alpha}{\sigma\over \sqrt{n}}-b$, then  $$L_{z2b}(\bar{x})= -\infty\,\ \mbox{and}\,\ U_{z2b}(\bar{x}) = \bar{x}-{\sigma\over \sqrt{n}} F^{-1}_N(\alpha-1+F_N(-{\sqrt{n}\over \sigma}(b+\bar{x}))). 
$$

The above discussion shows that the h-function method always generates an exact interval of level $1-\alpha$; however, depending on whether an appropriate point estimator is picked, the resultant interval can be the shortest z-interval or the widest interval, $(-\infty, +\infty)$.
~\raisebox{.5ex}{\fbox{}}

\medskip

\section{ Refining a $1-\alpha$ exact confidence interval}
When $C_0(\underline{X})$ is a $1-\alpha$ exact interval, what is the relationship between $C_0(\underline{X})$ and $C_0^M(\underline{X})$? From Tables~\ref{table-1} and \ref{table-2}, we see that the TIL decreases when the modification process is applied to an exact interval. Table~\ref{table-3} even shows that the modification process shrinks the interval uniformly on all sample points. So,
we show next that $C_0^M(\underline{X})$ is a subset of $C_0(\underline{X})$.

\subsection{Refining a $1-\alpha$ exact interval one time}

\begin{thm}
	\label{thm-ci-refine}
	For a $1-\alpha$ exact confidence interval $C_0(\underline{X})$, 
	interval
	$C_0^M(\underline{x})$ given in (\ref{test-ci-2}) is contained in	
	$C_0(\underline{X})$. Therefore, for any level confidence interval $C_0(\underline{X})$, $C_0^{M2}(\underline{X})$ is a subset of  $C_0^M(\underline{X})$.
\end{thm}

{\bf Example 1 (continued)}. It is clear in Table~\ref{table-3} that the modification of any 95\% exact interval is always a subset of the interval. We also apply Theorem~\ref{thm-ci-refine} to the four exact intervals in Table~\ref{table-1}: Clopper-Pearson interval $C_{p1}$ (1934),  Blaker interval (2000) $C_{p2}$, the
  exact likelihood-ratio-test interval $C_{p3}$ and Wang interval (2014) $C_{p4}$.  Table~\ref{table-4} contains the details. The improvements, as a proper subset, over $C_{p1}$ and $C_{p3}$ are noticeable. The improvement over $C_{p2}$ occurs at $x=6,8,10$. 
No improvement over $C_{p4}$ can be found since it is admissible for any $n$.  
In this sense, $C_{p4}$ is the best among the exact intervals discussed here.  In Tables~\ref{table-3} and \ref{table-4}, we provide five different admissible intervals for $p$,  $C_{p1}^{M}$,  $C_{p2}^{M}$,  $C_{p4}$, $C_{p6}^{M}$, 
$C_{p7}^{M}(=C_{p8}^{M17})$ with the same total interval length 6.4978. Four of them are generated by the newly proposed modification process.
~\raisebox{.5ex}{\fbox{}}

{\bf Example 3 (continued)}. Consider a variant of the z-interval for $\mu$: $C_0(\bar{X})=[\bar{X}-a{\sigma\over \sqrt{n}},\bar{X}+b{\sigma\over \sqrt{n}}]$ for two constants $a,b\geq z_{{\alpha\over 2}}$, which is a $1-\alpha$ exact but conservative interval. We now compute $C_0^M(\bar{X})$. 
Let $$T_0(\bar{X},\mu_0)=\min\{\mu_0-\bar{X}+ a{\sigma\over \sqrt{n}},\bar{X}+ b{\sigma\over \sqrt{n}}-\mu_0\}
.$$
\begin{table}[H]
	\begin{center}
		\renewcommand{\arraystretch}{0.9}
		\caption{ Clopper-Pearson interval $C_{p1}$, the refinement $C^M_{p1}$;  Blaker interval $C_{p2}$, the refinement $C^M_{p2}$; the likelihood-ratio-test interval $C_{p3}$, the refinement $C_{p3}^M$, the 22nd refinement $C_{p3}^{M22}$;     Wang interval $C_{p4}$; and their ICPs and TILs when $n=16$. For each interval, the first line is the lower limits and the second line is the upper limits. 
		} 
		\resizebox{1\textwidth}{0.35\textheight}{
			\begin{tabular}{l|rrrrrrrrrr}
				$X$  &   0    & 1      & 2     & 3      & 4      & 5      & 6      & 7      & 8 \\ \hline
				$C_{p1}$ & 0.0000 & 0.0015 & 0.0155& 0.0404 & 0.0726 & 0.1101 & 0.1519 & 0.1975 & 0.2465\\ 
				ICP=95.78\%& 0.2060 & 0.3024 & 0.3835& 0.4565 & 0.5238 & 0.5867 & 0.6457 & 0.7013 & 0.7535\\ 
				$C_{p1}^M(=C_{p1}^{M\infty})$& 0.0000 & 0.0032 & 0.0226 & 0.0531 & 0.0902 & 0.1321 & 0.1777 & 0.2017 & 0.2719\\
				ICP=95\% & 0.2018 & 0.3006 & 0.3690 & 0.4350 & 0.5000 & 0.5651 & 0.6311 & 0.6995 & 0.7281\\ \hline 
				$C_{p2}$ & 0.0000 & 0.0032 & 0.0226 & 0.0531 & 0.0902 & 0.1321 & 0.1746 & 0.2011 & 0.2717\\
				ICP=95\% & 0.2012 & 0.3005 & 0.3683 & 0.4345 & 0.5000 & 0.5656 & 0.6318 & 0.6996 & 0.7283\\ 
				$C_{p2}^M=(C_{p2}^{M\infty})$& 0.0000 & 0.0032 & 0.0226 & 0.0531 & 0.0902 & 0.1321 & 0.1777 & 0.2011 & 0.2719\\
				ICP=95\% & 0.2012 & 0.3005 & 0.3683 & 0.4345 & 0.5000 & 0.5656 & 0.6318 & 0.6996 & 0.7281\\ \hline
				$C_{p3}$ & 0.0000 & 0.0032 & 0.0226 & 0.0531 & 0.0902 & 0.1205 & 0.1462 & 0.1727 & 0.2592\\
				ICP=95\% & 0.1738 & 0.2885 & 0.3614 & 0.4312 & 0.5000 & 0.5689 & 0.6387 & 0.7116 & 0.7408\\ 
				$C_{p3}^M$& 0.0000 & 0.0032 & 0.0226 & 0.0531 & 0.0902 & 0.1321 & 0.1600 & 0.1732 & 0.2719\\
				ICP=95\% & 0.1738 & 0.2885 & 0.3614 & 0.4312 & 0.5000 & 0.5689 & 0.6387 & 0.7116 & 0.7281\\
				$C_{p3}^{M22}(=C_{p3}^{M\infty})$& 0.0000& 0.0032& 0.0226& 0.0531& 0.0902& 0.1321& 0.1737& 0.1737& 0.2719\\
				ICP=95\% & 0.1738& 0.2885& 0.3614& 0.4312& 0.5000& 0.5689& 0.6387& 0.7116& 0.7281\\
				\hline
				$C_{p4}(=C_{p4}^{M\infty})$ & 0.0000 & 0.0032 & 0.0226 & 0.0531 & 0.0902 & 0.1321 & 0.1777 & 0.2059 & 0.2719\\
				ICP=95\% & 0.2060 & 0.3024 & 0.3835 & 0.4416 & 0.5000 & 0.5585 & 0.6166 & 0.6977 & 0.7281\\ \hline
				$X$ &   9 & 10 & 11 & 12 & 13 & 14 & 15 & 16 & \\ \hline
				$C_{p1}$ & 0.2987 & 0.3543 & 0.4133 & 0.4762 & 0.5435 & 0.6165 & 0.6976 & 0.7940 \\ 
				TIL=6.9380 & 0.8025 & 0.8481 & 0.8899 & 0.9274 & 0.9596 & 0.9845 & 0.9985 & 1.0000\\ 
				$C_{p1}^M$& 0.3005 & 0.3689 & 0.4349 & 0.5000 & 0.5650 & 0.6310 & 0.6994 & 0.7982\\ 
				TIL=6.4978 & 0.7983 & 0.8223 & 0.8679 & 0.9098 & 0.9469 & 0.9774 & 0.9968 & 1.0000\\ \hline
				$C_{p2}$   & 0.3004 & 0.3682 & 0.4344 & 0.5000 & 0.5655 & 0.6317 & 0.6995 & 0.7988\\
				TIL=6.5043 & 0.7989 & 0.8254 & 0.8679 & 0.9098 & 0.9469 & 0.9774 & 0.9968 & 1.0000\\ 
				$C_{p2}^M$ & 0.3004 & 0.3682 & 0.4344 & 0.5000 & 0.5655 & 0.6317 & 0.6995 & 0.7988\\
				TIL=6.4978 & 0.7989& 0.8223 & 0.8679 & 0.9098 & 0.9469 & 0.9774 & 0.9968 & 1.0000\\ \hline
				$C_{p3}$   & 0.2884 & 0.3613 & 0.4311 & 0.5000 & 0.5688 & 0.6386 & 0.7115 & 0.8262\\
				TIL=6.6115 & 0.8273 & 0.8538 & 0.8795 & 0.9098 & 0.9469 & 0.9774 & 0.9968 & 1.0000\\ 
				$C_{p3}^M$ &0.2884 & 0.3613 & 0.4311 & 0.5000 & 0.5688 & 0.6386 & 0.7115 & 0.8262\\
				TIL=6.5343 & 0.8268 & 0.8400 & 0.8679 & 0.9098 & 0.9469 & 0.9774 & 0.9968 & 1.0000\\
				$C_{p3}^{M22}$& 0.2884& 0.3613& 0.4311& 0.5000& 0.5688& 0.6386& 0.7115& 0.8262\\
				TIL=6.5058 & 0.8263& 0.8263& 0.8679& 0.9098& 0.9469& 0.9774& 0.9968& 1.0000\\ \hline
				$C_{p4}$   & 0.3023 & 0.3834 & 0.4415 & 0.5000 & 0.5584 & 0.6165 & 0.6976 & 0.7940\\
				TIL=6.4978 & 0.7941& 0.8223 & 0.8679 & 0.9098 & 0.9469 & 0.9774 & 0.9968 & 1.0000\\
			\end{tabular}
		}
		\label{table-4}
	\end{center}
\end{table}

For a fixed $\bar{x}$, $T_0(\bar{x},\mu_0)$ is unimodal in $\mu_0$ with the mode at $\mu_0=\bar{x}+{\sigma\over \sqrt{n}}({b-a\over2 })$. 
Algebraic calculation shows, for $Z\sim N(0,1)$, 
\begin{eqnarray*}
	&&h_0(\bar{x},\mu_0)
	= P(\min\{-Z+a,Z+b\}\leq {\sqrt{n}\over \sigma}T_0(\bar{x},\mu_0))\\
	&=&\left\{
	\begin{split}
		& P(Z>\max\{{a-b\over 2}, - {\sqrt{n}\over \sigma}(\mu_0-\bar{x})\})+P(Z<\min\{{a-b\over 2}, \\ & \hspace{1in} a-b+{\sqrt{n}\over \sigma}(\mu_0-\bar{x})\})	 \,\ \mbox{if} \,\  \mu_0\leq \bar{x}+{\sigma\over \sqrt{n}}({b-a\over2 }),\\
		&  P(Z>\max\{{a-b\over 2}, {\sqrt{n}\over \sigma}(\mu_0-\bar{x})+a-b\})+P(Z<\min\{{a-b\over 2},\\ & \hspace{1in} -{\sqrt{n}\over \sigma}(\mu_0-\bar{x})\})	 \,\ \mbox{if} \,\  \mu_0> \bar{x}+{\sigma\over \sqrt{n}}({b-a\over2 }).\\ 
	\end{split}
	\right.
\end{eqnarray*}
This function is also unimodal in $\mu_0$ with the same mode of $T_0(\bar{x},\mu_0)$. Note that $h_0$ at the model is equal to one and goes to zero as $\mu_0$ goes to infinity. Let $L_0^M(\bar{x})$ and $U_0^M(\bar{x})$ be the two solutions of $h_0(\bar{x},\mu_0)=\alpha$. The two solutions can be written as $\bar{x}+c_i(a,b){\sigma\over \sqrt{n}}$ for $i=1$ and 2, where $c_1(a,b)=-z_{\alpha_1}$ and $c_2(a,b)=z_{\alpha_2}$ for  $\alpha_1+\alpha_2=\alpha.$ When $a=b$, then $C_0^M$ is equal to the z-interval.
~\raisebox{.5ex}{\fbox{}}

\subsection{Refining a $1-\alpha$ exact interval multiple times}
The modification process can be applied to a $1-\alpha$ exact interval repeatedly, and in each time an interval which is a subset of the previous interval is generated. 
Next we determine the smallest $1-\alpha$ exact interval in the modification process.
\begin{thm}
	\label{thm-ci-refine-inf}
	For a $1-\alpha$ exact confidence interval $C_0(\underline{X})$ for $\theta$, which is rewritten as $C_0^{Mk}(\underline{X})$ for $k=0$, let $C_0^{M\infty}(\underline{x})=\cap_{k=0}^{+\infty}C_0^{Mk}(\underline{x})$ for any $\underline{x}$.  Then,
	
	i) interval
	$C_0^{Mk}(\underline{x})$, as a set of $\theta$, is  nonincreasing in $k$ for $k\geq 0$.

	ii) $C_0^{M\infty}(\underline{X})$, contained in $C_0^{Mk}(\underline{X})$ for any $k$, is  a $1-\alpha$ exact confidence interval. 
	
	iii) If $C_0^{Mk}(\underline{X})=C_0(\underline{X})^{M(k+1)}$ for some $k\geq 0$, then 
	$C_0^{M\infty}(\underline{X})=C_0^{Mk}(\underline{X})$.
\end{thm}

When applying Theorem~\ref{thm-ci-refine-inf}, the smallest interval $C^{M\infty}(\underline{X})$ is typically  equal to $C^{Mk}(\underline{X})$ for a finite nonnegative integer $k$.   For examples, in Table~\ref{table-3}, $C_{p6}^{M\infty}=C_{p6}^{M}$ for $k=1$, and in Table~\ref{table-4}, $C_{p4}^{M\infty}=C_{p4}$ for $k=0$.

\medskip
{\bf Example 2 (continued)}. 
Following a series of works, originated by Buehler (1957), the construction of the smallest $1-\alpha$ one-sided confidence interval is an automatic process as long as an order on the sample space is prespecified. Therefore, a $1-\alpha$ exact two-sided confidence interval can be easily obtained by taking the intersection of two lower and upper $1-{\alpha\over 2}$ one-sided intervals. When the sample space is discrete, the resultant interval may be conservative with an ICP higher than the nominal level $1-\alpha$ (Agresti, 2013).  Clopper-Pearson interval (1934) is such an example. There are successful efforts to improve this interval uniformly, see Casella (1986), Blaker (2000) and Wang (2014), since a single binomial distribution  does not involve any nuisance parameter. When there exists a nuisance parameter, as in many applications including the inference about the difference $d$, it is challenging to improve the intersection of two $1-{\alpha\over 2}$ one-sided intervals.
The main difficulty is that the method of finding an ICP in Wang (2007) fails.  Theorem~\ref{thm-ci-refine}, however, provides a promising effort to shrink any exact interval $C_0(\underline{X})$, including the intersection of two $1-{\alpha\over 2}$ one-sided intervals. 

To be precise, $C_{d5}(X,Y)$ in Table \ref{table-2} is the  intersection of two $1-{\alpha\over 2}$ smallest one-sided intervals in Wang (2010). Its ICP and TIL are larger than those of $C_{d2}(X,Y)$, but $C_{d5}^{M\infty}(X,Y)$ has the smallest TIL in general. The reason is that the confidence limits of $C_{d5}(X,Y)$ do not have many ties due to an inductive order on the sample space. On the contrary, $C_{d4}(X,Y)$, which is equal to the MLE estimator $\hat{p}_1-\hat{p}_2$, has many ties. For example, when $(n_1,n_2)=(8,10)$, then 
$C_{d4}(0,0)=C_{d4}(4,5)=C_{d4}(8,10)=0$. This results in the tied intervals of $C^{M\infty}_{d4}(0,0)=C^{M\infty}_{d4}(4,5)=C^{M\infty}_{d4}(8,10)=[-0.4375, 0.4375]$. In consequence, $C^{M\infty}_{d4}(X,Y)$ is much longer than
$C_{d5}^{M\infty}(X, Y)$. Score interval $C_{d2}(X,Y)$ also has a small TIL due to a small number of ties. 

How to determine whether $C_{0}^{Mk}$ for an integer $k$ is equal to $C_{0}^{M\infty}$?
Since the intervals for $d$ all have a finite length, the ratio of the TILs of two consecutive $C_{0}^{M(k+1)}$ and $C_{0}^{Mk}$ is used. Due to part i) of Theorem~\ref{thm-ci-refine}, the ratio is less than or equal to one. If it is equal to one, then, by part iii) of Theorem~\ref{thm-ci-refine}, $C_{0}^{Mk}=C_{0}^{M\infty}$. 
For example, $C_{d2}^{M\infty}=C_{d2}^{M20}$ and $C_{d5}^{M\infty}=C_{d5}^{M17}$. Here the ratio is accurate up to the seventh decimal place.
~\raisebox{.5ex}{\fbox{}}

\medskip
{\bf Example 4}. Consider a two-arm randomized trial in Essenberg (1952) for testing
the effect of tobacco smoking on tumor development in mice.
The smoking group had $23(= n_1)$ mice and tumors were observed on $21(= x)$ mice; $(n_2,y) = (32,19)$  in the control group. Let $p_1$ and $p_2$ be the tumor rates
for the two groups. We use $d$ to evaluate the smoking effect.
The fifteen intervals in Table~\ref{table-2} at $(x,y)=(21,19)$ are reported in 
Table~\ref{table-5}. The final refinement of Wang interval, $C^{M\infty}_{d5}(x,y)$, is equal to $[0.0968, 0.5037]$ and has the shortest length 0.4049, an ICP 0.95 and the shortest TIL 341.4697. ~\raisebox{.5ex}{\fbox{}}
\begin{table}[h]
	\begin{center}
		\renewcommand{\arraystretch}{0.9}
		\caption{Fifteen 95\% confidence intervals at $(x,y)=(21,19)$: $(C_{di}$, $C^M_{di}$, $C^{M\infty}_{di})$ for $i=1$ (LRT, exact), $2$ (Score, exact), $3$ (Wald, approximate), $4$ (MLE as an interval) and $5$ (Wang, exact) and their lengths when $(n_1,n_2)=(23,32)$. The smallest length of the exact intervals is marked by `*'.} 
		\resizebox{1\textwidth}{0.1\textheight}{
			\begin{tabular}{cccc|cccc|cccc}
				& lower & upper & length & & lower & upper & length &  & lower & upper & length  \\ \hline
				$C_{d1}$ &  0.0607 & 0.5337 & 0.4730 & $C^M_{d1}$ & 0.0610 & 0.5337 & 0.4727 & $C^{M\infty}_{d1}$ & 0.0612 & 0.5337  & 0.4725\\
				$C_{d2}$ &  0.0794 & 0.5227 & 0.4433 &$C^M_{d2}$ & 0.0794 & 0.5222 & 0.4428 & $C^{M\infty}_{d2}$ & 0.0794 & 0.5217 & 0.4423\\
				$C_{d3}$ &  0.1138 & 0.5248 & 0.4110 &$C^M_{d3}$ & 0.0569 & 0.5485 & 0.4916 & $C^{M\infty}_{d3}$ & 0.1140 & 0.5470 & 0.4330\\
				$C_{d4}$ &  0.3193 & 0.3193 & 0      & $C^M_{d4}$ & 0.0523& 0.5442 & 0.4919 & $C^{M\infty}_{d4}$ & 0.0530 & 0.5438  & 0.4908\\
				$C_{d5}$ &  0.0947 & 0.5126 & 0.4179 & $C^M_{d5}$ & 0.0968& 0.5080 & 0.4112 & $C^{M\infty}_{d5}$ & 0.0968 & 0.5037  & 0.4069*\\
			\end{tabular}
			\label{table-5}
		}
	\end{center}
\end{table}

\subsection{Confidence intervals for the difference of two proportions in a matched-pair experiment} Consider a 2$\times$2 contingency table with two binary variables $A$ (row) and $B$ (column), where $1$ is a success and $0$ is a failure. 
We observe a random vector $(N_{11},N_{10},N_{01},N_{00})$ that follows a multinomial distribution with $n$ trials and probabilities $(p_{11},p_{10},p_{01},p_{00})$.
For example, $N_{10}$ is the number of subjects with $(A,B)=(1,0)$ in all $n$ subjects in a study and $p_{10}=P((A,B)=(1,0))$. Similar to $d$ in Section 2.3, the parameter of interest here is $$d_m=P(A=1)-P(B=1)=p_{10}-p_{01}.$$ Let $T=N_{11}+N_{00}$ and $p_t=p_{11}+p_{00}$. The conditional distribution of $N_{ij}$'s for given $(N_{10},T)$
does not involve $p_{10}$ and $p_{01}$, so the inferences about $d_m$ should be based on $(N_{10},T)$ following the similar reasoning of the sufficiency principle. 
The reduced sample and parameter spaces are $S_M=\{(n_{10},t): n_{10}+t\in[0,n] \}$ and $H_M=\{(d_m,p_t):  d_m\in [-1,1], p_t \in [0,1-|d_m|]\},$
respectively.
The PMF for $(N_{10},T)$, in terms of $(d_m,p_t)$, is
$$p_{M}(n_{10},t,d_m,p_t)={n! \over n_{10}!t!n_{01}! } ({1+d_m-p_t\over 2})^{n_{10}}p_t^t({1-d_m-p_t\over 2})^{t-n_{10}}.$$ 

Wang (2012) proposed the smallest  $1-\alpha$ lower and upper one-sided intervals for $d_m$ using an inductive order on $S_M$. Then, the intersection of two smallest  $1-{\alpha \over 2}$ lower and upper one-sided intervals for $d_m$, denoted by $C_{d_m1}(N_{10},T)=[L_{d_m1}(N_{10},T),U_{d_m1}(N_{10},T)]$, is  of level $1-\alpha$. This interval can be computed by an R-package ``ExactCIdiff'' (Shan and Wang, 2013). It is now refined following the modification process. More precisely,
let 
\begin{equation}
\label{ts-dm1}
T_{m1}(N_{10},T)=\min\{d_m-L_{d_m1}(N_{10},T),U_{d_m1}(N_{10},T)-d_m\}.
\end{equation}
Define an h-function
\begin{equation}\label{h-dm1}
h_{m1}(n_{10},t,d_m)
= \sup_{p_t\in [0,1-|d_m|]} \sum_{\{(n_{10}',t')\in S_M: T_{m1}(n_{10}',t')\leq T_{m1}(n_{10},t)\}} p_{M}(n_{10},t,d_m,p_t).
\end{equation}
Then, the interval generated by the modification process 
\begin{equation}\label{ci-dm1}
C_{d_m1}^M(n_{10},t)=\overline{\{d_m\in [-1,1]: h_{m1}(n_{10},t,d_m)>\alpha\}}
\end{equation}
is of level $1-\alpha$ and is a subset of $C_{d_m1}(n_{10},t)$ following Theorem~\ref{thm-ci-refine}. We continue the modification process for $k$ times so that no improvement can be found. i.e,
$C_{d_m1}^{Mk}=C_{d_m1}^{M\infty}$. 

There have been many efforts to derive approximate and exact confidence intervals for $d_m$.  Fagerland, Lydersen and Laake (2014) provided a good summary. Here are two main results: i) Tango (1998) recommended using  his approximate score interval for large
$|d_m|$ and two intervals in Newcombe (1998) that are based on the Wilson score interval or Jeffreys interval for the single proportion for small $|d_m|$. The formula of the first interval, denoted by $C_{d_m2}$, is given in Fagerland et al. (2014, equations 24, 25). The other two intervals depend on 
$(N_{11},N_{10},N_{01})$ but not $(N_{10},T)$. So, any of the two intervals assumes different interval-values on those sample points of $(n_{11},n_{10},n_{01})$ with a fixed $(n_{10},t)$. In consequence, this may cause an interpretation problem of the intervals. The two intervals are not included in the comparison below.  ii) Bonett and Price (2012) compared five intervals and concluded that the Wald interval with Bonett–Price adjustment, denoted by $C_{d_m3}$ (Fagerland et al. 2014, equation 16), performs as well as or better than $C_{d_m2}$. These two intervals are not of level $1-\alpha$. However, 
following the similar process of (\ref{ts-dm1}), (\ref{h-dm1}) and (\ref{ci-dm1}) for deriving 
$C_{d_m1}^{M}$ and $C_{d_m1}^{M\infty}$, the counterpart intervals for $C_{d_m2}$ and $C_{d_m3}$ are constructed. They are $1-\alpha$ exact intervals due to Theorems~\ref{thm-ci-refine} and \ref{thm-ci-refine-inf}.
Next, we present a limited numerical comparative study for $C_{d_m1}$, $C_{d_m2}$, $C_{d_m3}$ and their modifications in terms of ICP and TIL. 
\medskip

{\bf Example 5}. Bentur et al. (2009) measured airway hyper-responsiveness status: (Yes = 1, No = 0)  in $n(=21)$ children before (A) and after (B) stem cell transplantation and observed $(n_{11},n_{10},n_{01},n_{00})=(1,1,7,12)$ and $t=13$. Then, the MLE estimate for $d_m$ is $-0.2857$.  In Table~\ref{table-6}, we report nine 95\% confidence intervals at the observation for an individual performance and their ICPs and TILs for an overall performance. 

As expected, $C_{d_m1}$ has the largest length and TIL. 
However, the small TILs for $C_{d_m2}$ and $C_{d_m3}$ are due to their incorrect ICPs: 0.8376 and 0.9146. After the modification process is applied to the three intervals, the resultant intervals all have ICPs no less than 0.95; $C_{d_m1}^{M\infty}$ is the shortest at the observation and has a little larger TIL than $C_{d_m2}^{M\infty}$. One reason for the large TIL of $C_{d_m3}^{M\infty}$ is that $C_{d_m3}$ has ties in its confidence limits, especially when $n_{10}$ is close to $n$ or $0$.  For examples, the upper limits for $C_{d_m3}$ at $(n_{10},t)=(21,0), (20,0), (20,1), (19,0), (19,1)$ and (19,2) are equal to 1. In order to obtain the final interval, the numbers of modification, $k$, for the three intervals are equal to 18, 19, and 20, respectively. When $k$ goes large, the improvement in TIL is negligible. For example, when $k=5$, $C_{d_m2}^{Mk}$ has a TIL of 146.2981, close to its final TIL 146.2317.  
\begin{table}[h]
	\begin{center}
		\caption{Nine 95\% confidence intervals: $(C_{d_mi}$, $C^M_{d_mi}$, $C^{M\infty}_{d_mi})$ for $i=1$ (Wang, exact), $2$ (Tango, approximate), $3$ (Wald with Bonett and Price adjustment, approximate), their lengths at $(n_{10},t)=(1,13)$, their ICPs and TILs when $n=21$. The smallest length and TIL of the exact intervals are marked by `*'.} 
		\resizebox{1\textwidth}{0.25\textheight}{
			\begin{tabular}{cccc|cccc|cccc}
				\multicolumn{12}{c}{Part I: The nine intervals 
					at $(n_{10},t)=(1,13)$}\\
				& lower & upper & length & & lower & upper & length &  & lower & upper & length  \\
				& \multicolumn{2}{c}{} && &\multicolumn{2}{c}{the p-value} & & &\multicolumn{2}{c}{the p-value} &\\
				\hline
				$C_{d_m1}$ &  -0.5214& -0.0126 & 0.5088& $C^M_{d_m1}$ & -0.5065 & -0.0155 & 0.4910 & $C^{M\infty}_{d_m1}$ & -0.4923& -0.0155& 0.4768* \\
				& \multicolumn{2}{c}{} && &\multicolumn{2}{c}{0.04125} && &\multicolumn{2}{c}{0.04393} &\\
				$C_{d_m2}$ & -0.5173& -0.0260 & 0.4913 &$C^M_{d_m2}$ & -0.5320& -0.0182 & 0.5138 & $C^{M\infty}_{d_m2}$ & -0.5287 & -0.0182& 0.5105 \\
				& \multicolumn{2}{c}{} & &&\multicolumn{2}{c}{0.04125} && &\multicolumn{2}{c}{0.04215} &\\				
				$C_{d_m3}$ & -0.5084& -0.0133 & 0.4951 &$C^M_{d_m3}$ & -0.5000 &  0.0122&  0.5122  & $C^{M\infty}_{d_m3}$ & -0.4997 & 0.0122 & 0.5119 \\ 
				& \multicolumn{2}{c}{} & &&\multicolumn{2}{c}{0.07835} && &\multicolumn{2}{c}{0.07835} &\\ 
				\\
				\multicolumn{12}{c}{Part II: The nine interval's ICPs and TILs}\\
				& ICP & TIL &  & & ICP & TIL &  &  & ICP & TIL &   \\ \hline
				$C_{d_m1}$ &  0.9500 & 152.4780 &  & $C^M_{d_m1}$ & 0.9500 & 147.8739 &  & $C^{M\infty}_{d_m1}$ & 0.9500 & 146.8296  & \\
				$C_{d_m2}$ &  0.8376 & 144.1614 &  &$C^M_{d_m2}$ & 0.9500 & 147.7267 &  & $C^{M\infty}_{d_m2}$ & 0.9500 & 146.2317* & \\
				$C_{d_m3}$ &  0.9146 & 147.7201 &  &$C^M_{d_m3}$ & 0.9501 & 152.3374 & & $C^{M\infty}_{d_m3}$ & 0.9500 & 149.2308 &\\ \\
			\end{tabular}
			\label{table-6}
		}
	\end{center}
\end{table}

In Table~\ref{table-6}, we also report the p-values for testing $H_0: d_m=0$ that are associated with the modified intervals $C_{d_mi}^M$ and $C_{d_mi}^{M\infty}$ at $(n_{10},t)=(1,13)$ for $i=1,2,3$. For example, the p-value corresponding to $C_{d_m1}^M$ is equal to $h_{m1}(1,13,0)=0.04125$, where $h_{m1}$ is given in (\ref{h-dm1}). We reject $H_0$ at level 0.05, which is consistent to the fact that interval $C_{d_m1}^M(1,13)$ excludes zero; while $C_{d_m3}^M(1,13)$ includes zero and fails to reject $H_0$ with a p-value of 0.07835 although $C_{d_m3}(1,13)$ excludes zero. 
~\raisebox{.5ex}{\fbox{}}

\section{Modifying one-sided intervals}

Assume the range of parameter $\theta$ is $[A,B]$ for two known constants $A$ and $B$. If one end is infinity, for example, $B=+\infty$, then $[A,B]=[A,+\infty)$.
For a lower one-sided confidence interval for $\theta$, $C_l(\underline{X})=[L_l(\underline{X}),B]$ with $L_l(\underline{X})\geq A$, 
consider the one-sided hypotheses 
\begin{equation}\label{hy-l}
H_0:\theta\leq \theta_0\,\ vs \,\ H_A:\theta>\theta_0.
\end{equation}
Let $$T_{1l}(\underline{x},\theta_0)=\theta_0-L_l(\underline{x}).$$  
A small value of $T_{1l}$ gives strong evidence to establish $H_A$,
and $$T_{1l}(\underline{x},\theta_0)\geq 0 \,\ \mbox{if and only if} \,\ \theta_0 \geq L_l(\underline{x})\,\ (i.e., \,\ \theta_0\in C_l(\underline{x})).$$ 
The h-function based on $T_{1l}(\underline{x},\theta_0)$ is
\begin{equation}
\label{h-1l}
h_{1l}(\underline{x},\theta_0)=\sup_{H_0}P(T_{1l}(\underline{X},\theta_0)\leq T_{1l}(\underline{x},\theta_0))=\sup_{\theta\leq \theta_0}P(L_l(\underline{x})\leq L_l(\underline{X})).
\end{equation}
Following (\ref{test-ci}), the level-$\alpha$ acceptance region for $H_0$ and the  $1-\alpha$ exact lower one-sided confidence interval for $\theta$ are 
\begin{equation}
\label{test-ci-diff-1l}
A_{1l}(\theta_0)=\{\underline{x}: h_{1l}(\underline{x},\theta_0)>\alpha\}\,\ \mbox{and}\,\ C_l^M(\underline{X})=\overline{\{\theta_0: h_{1l}(\underline{x},\theta_0)>\alpha\}}.
\end{equation}
As mentioned in Example 2 (continued), if an order is introduced on the sample space, then the smallest $1-\alpha$ one-sided confidence interval under this order can be automatically constructed. 
The order determines, between any two sample points, which point has a larger confidence limit and can be given by a function on the sample space. In the current case, the function is equal to $L_l(\underline{X})$. More precisely, consider a $1-\alpha$ exact lower one-sided interval class ${\cal C}_l=\{C(\underline{X})=[L(\underline{X}),B]: L(\underline{x}')\leq(=) L(\underline{x})\,\ \mbox{ if}\,\ L_l(\underline{x}')\leq (=)L_l(\underline{x}), \,\ \forall\,\ \underline{x}' \,\ \mbox{and} \,\ \underline{x} \}$. 
We should find the smallest interval in ${\cal C}_l$, which is contained in any interval in ${\cal C}_l$.
Different from the two-sided interval case, we have the following stronger result than
Theorems~\ref{thm-ci-modi}, \ref{thm-ci-refine} and \ref{thm-ci-refine-inf}. One time of the modification yields the smallest interval. This also establishes a connection between the h-function method and the construction of the smallest one-sided interval under an order.

\begin{thm}
	\label{thm-ci-1l}
	For a lower one-sided confidence interval $C_l(\underline{X})=[L_l(\underline{X}), B]$ of any level, 
	
	i) interval
	$C_l^M(\underline{X})$ given in (\ref{test-ci-diff-1l}) is of level $1-\alpha$;
	
	ii) $C_l^{M\infty}(\underline{X})=C_l^M(\underline{X})$; 
	
	iii)  $C_l^M(\underline{X})=[L_l^M(\underline{X}),B]$ is the smallest interval in ${\cal C}_l$. i.e., for any interval $C(\underline{X})$ in ${\cal C}_l$, $L_l^M(\underline{X})\geq L(\underline{X})$. 
\end{thm}

{\bf Example 7}. (the stochastically nondecreasing distribution family with a single parameter)
Suppose $X$ has a CDF $F(x,\theta)$ and is stochastically nondecreasing in $\theta$. i.e., 
\begin{equation}\label{sto-greater}
F(x,\theta_1)\geq F(x,\theta_2)\,\ \mbox{ for any $x$ and $\theta_1\leq \theta_2$}. 
\end{equation}
This family includes all important single-parameter distributions, such as binomial, hypergeometric, exponential, Poisson, etc. 
We are interested in finding the modified interval $C_l^M(X)$ for the one-sided confidence interval $C_l(X)=[X,B]$. Interval $C_l(X)$ itself may be meaningless in terms of  estimating $\theta$. For example, when 
$X\sim Bino(n,p)$, no one would use $[X,1]$ to estimate $p$. However, a large value of $X$ 
tends to associate with a large value of $\theta$ as seen in the monotone likelihood ratio distribution family, a special case of the stochastically nondecreasing distribution family.
Therefore, the observation $X$ is able to tell which sample point has a large lower limit. This is the only fact needed to derive  $C_l^M(X)$. 

To be precise, following (\ref{h-1l}), $$h_{1l}(x, \theta_0)=\max_{\theta\leq \theta_0}P(x\leq X)\stackrel{due\,\ to \,\ (\ref{sto-greater})}{=}1-F(x-1,\theta_0).$$
Note that $F(x,\theta)$ is nonincreasing in $\theta$ for any $x$. Then, $L_l^M(x)=\inf\{\theta_0: 1-F(x-1,\theta_0)>\alpha \}.$
This interval is the smallest interval among $1-\alpha$ exact intervals of form $[L(X),B]$ with a nondecreasing $L$.

For example, when $X$ follows $Bino(n,p)$, then $L_l^M(x)=0$ for $x=0$; and $L_l^M(x)$ is the solution of
$ 1-F_B(x-1,n,p_0)=\alpha$ for $x>0$. In fact,  $C_l^M(X)=[L_l^M(X),1]$ is equal to the $1-\alpha$ lower one-sided Clopper-Pearson interval. 
~\raisebox{.5ex}{\fbox{}}
\medskip

{\bf Example 8} (the one-sided t-interval). When estimating a normal mean $\mu$ with an unknown variance $\sigma^2$, consider a lower one-sided interval, $C_c(\bar{X},S)=[\bar{X}+c {S\over\sqrt{n}},+\infty)$ for a given constant $c$ and the sample standard deviation $S$. This is equal to the one-sided t-interval if $c=-t_{\alpha,n-1}$, where $t_{\alpha,n-1}$ is the upper $\alpha$th percentile a t-distribution with $n-1$ degrees of freedom. Now we derive $C_c^M(\bar{X},S)$. Note 
$\bar{X}\sim \mu+{\sigma\over \sqrt{n}}Z$ and $S^2 \sim {\sigma^2\over n-1}\chi^2_{n-1}$.
Following (\ref{h-1l}),
\begin{eqnarray*}
	h_{1l}(\bar{x},s,\mu_0)
	&=&\max_{\sigma>0}P({Z \over \chi_{n-1}/\sqrt{n-1}}  \geq {\sqrt{n}(\bar{x}+c{s\over\sqrt{n}}-\mu_0)\over \sigma \chi_{n-1}/\sqrt{n-1}}-c)\\
	&=& \left\{
	\begin{split}
		&1 \,\ \,\ \,\ \,\ \,\ \,\ \,\ \,\ \,\ \,\  \mbox{if}\,\ \mu_0>\bar{x}+c{s\over\sqrt{n}};\\
		&1-F_T(-c)\,\ \mbox{if}\,\ \mu_0\leq \bar{x}+c{s\over\sqrt{n}},
	\end{split}
	\right.
\end{eqnarray*} 
where $F_T$ is the CDF of t-distribution with $n-1$ degrees of freedom. Thus,
$$ C_c^M(\bar{x},s)=\overline{\{\mu_0: h_{1l}(\bar{x},s,\mu_0)>\alpha\}}=\left\{
\begin{split}
& (-\infty,+\infty) \,\ \mbox{if} \,\ c>-t_{\alpha,n-1};\\
& [\bar{x}+c{s\over\sqrt{n}},+\infty)  \,\ \mbox{if} \,\ c\leq -t_{\alpha,n-1},
\end{split}
\right.
$$  by (\ref{test-ci-diff-1l}). Furthermore, among all intervals $C_c^M$ for different $c$'s, the one with $c=-t_{\alpha,n-1}$ is the smallest one.
~\raisebox{.5ex}{\fbox{}}
\medskip

The above two examples show the importance of selecting a good order function on interval construction from both the positive and negative sides. In Example 7, $C_l$ may not be a meaningful interval, but the good order by $X$ still generates the smallest interval; in Example 8, the order by $\bar{X}+c{S\over \sqrt{n}}$ is not good, except the case of $c=-t_{\alpha,n-1}$, and does not yield an optimal interval. 

We next describe a parallel result
for an upper one-sided confidence interval $C_u(\underline{X})=[A, U(\underline{X})]$.
Let  $$T_{1u}(\underline{x},\theta_0)=U_u(\underline{x})-\theta_0.$$  
Then, the h-function based on $T_{1u}(\underline{x},\theta_0)$ is
\begin{equation}
\label{h-1u}
h_{1u}(\underline{x},\theta_0)=\sup_{H_0}P(T_{1u}(\underline{X},\theta_0)\leq T_{1u}(\underline{x},\theta_0))=
\sup_{\theta\geq \theta_0}P(U_l(\underline{X})\leq U_l(\underline{x})).
\end{equation}
Following (\ref{test-ci}), the level-$\alpha$ acceptance region for $H_0$ and  $1-\alpha$ exact upper one-sided confidence interval for $\theta$ are 
\begin{equation}
\label{test-ci-diff-1u}
A_{1u}(\theta_0)=\{\underline{x}: h_{1u}(\underline{x},\theta_0)>\alpha\}\,\ \mbox{and}\,\ C_u^M(\underline{x})=\overline{\{\theta_0: h_{1u}(\underline{x},\theta_0)>\alpha\}}.
\end{equation}
We have a similar result to Theorem~\ref{thm-ci-1l} and skip the proof.

\begin{thm}
	\label{thm-ci-1u}
	For an upper one-sided confidence interval $C_u(\underline{X})=[A, U_u(\underline{X})]$ of any level, 
	
	i) interval
	$C_u^M(\underline{X})$ given in (\ref{test-ci-diff-1u}) is of level $1-\alpha$;
	
	ii) $C_u^{M\infty}(\underline{X})=C_u^M(\underline{X})$. 
	
	iii) Define a $1-\alpha$ exact upper one-sided interval class ${\cal C}_u=\{C(\underline{X})=[A,U(\underline{X})]: U(\underline{x}')\leq U(\underline{x})\,\ \mbox{if}\,\ U_u(\underline{x}')\leq U_u(\underline{x}), \,\ \forall\,\ \underline{x}' \,\ \mbox{and} \,\ \underline{x} \}$. Then
	$C_u^M(\underline{X})=[A,U_u^M(\underline{X})]$ is the smallest interval in ${\cal C}_u$. i.e., for any interval $C(\underline{X})$ in ${\cal C}_u$, $U_u^M(\underline{X})\leq U(\underline{X})$. 
\end{thm}
\section{Discussions}

A confidence interval can be obtained by converting a family of tests and vice versa. 
In this paper, however, we formally introduce a middle function, the h-function, with which both the confidence interval and test can be derived. Thus, we propose a simpler but more general approach. This idea was used by researchers, such as Spjotvoll (1983), Blakers (2000) and Agresti (2013, p. 609), but the process is not clearly defined in the general setting as in the current paper. Casella and Berger (2002) introduced seven methods to construct tests and confidence intervals in great detail. The proposed h-function method is a supplement or replacement for some of those methods.

We show the effectiveness of the h-function method in Sections 2 through 5. The method is general and can be used in many applications, especially when a nuisance parameter exists.  Theorem~\ref{thm-ci-modi} is easy to follow and powerful in its ability to modify any confidence interval, including asymptotic intervals and point estimators, to exact intervals.  A straightforward extension of the theorem is to modify a credible interval to a $1-\alpha$ exact interval because the credible interval can be treated as a confidence interval. 
This greatly enhances the reliability of these intervals when they are applied to a finite sample since an invalid inferential procedure is converted to be a valid one. 
More importantly, when nuisance parameters are presented, it is challenging to improve an existing exact confidence interval. Theorem~\ref{thm-ci-refine} is a successful effort that delivers an exact interval which is a subset of the original interval, especially when the underlying distribution is discrete. Theorem~\ref{thm-ci-refine-inf} provides the smallest interval that can be generated by the modification process. Example 6 illustrates how the process modifies a conservative interval to a short one with an appropriate ICP. 
Further, from a mathematical point of view, we establish a connection between the h-function method and the construction of the smallest one-sided intervals based on an order presented within Theorems~\ref{thm-ci-1l} and \ref{thm-ci-1u}. 

From a theoretical point of view, the interval construction now becomes an automatic process. However, this process is still computationally complex,  particularly for precisely finding the globe maximum of a function, $P(K(\underline{x},\theta_0))$ in (\ref{h-function-ts-small}), as a function of $\underline{\eta}$, and solving the smallest and largest roots of the equation, $h(\underline{x},\theta_0)-\alpha>0$.   To our best knowledge, there is no effective computing software that can accomplish the two tasks both quickly and accurately. It is worth mentioning that $h(\underline{x},\theta_0)$ is not continuous in $\theta_0$ in general. Our best effort for global optimization is based on the combination of a grid search and a local optimization. This reduces down to a question of whether the resultant interval (e.g., $C_0^M$) is truly level $1-\alpha$ due to a grid search not being fine enough. We intend to select a large number of points for the search, which inevitably takes more computation time. For example, this number is between 200 and 2000 for a range of $[-1,1]$ when deriving intervals for $d$. Additionally, in each of the tables we also report the ICP to assure that each exact interval has a precise ICP.  Due to the theorems, large numbers in grid search and ICP confirmations, we are confident in the numerical results presented in the paper.          

As shown in Example 1, the two-sided interval $C_0^{M\infty}$ may or may not be admissible.  An interesting question is when  $C_0^{M\infty}$ is admissible. One can see $C_{p5}^{M\infty}$ in Table~\ref{table-3} is not admissible due to the tied values of $L_{p5}^{M\infty}(x)=0$ for $x=0,1,2,3$ and others. 
If $C_0^{M\infty}$ is not admissible, it cannot be improved by the proposed modification process. We may have two options: i)
Follow Casella (1986) or Wang (2014) to construct admissible intervals. However, 
the implementation is not clear if there exists a nuisance parameter. ii) Define an appropriate interval class so that  $C_0^{M\infty}$ is admissible within the class. Class ${\cal C}_l$ introduced before Theorem~\ref{thm-ci-1l}  is a successful example for the case of one-sided intervals.
We conjecture that if $C_0^{M\infty}$ has no ties then $C_0^{M\infty}$ is admissible.

Another research question is whether there exists a better choice than $T_2$ in (\ref{ts-t2}) so that $C_0^{Mk}$ converges to $C_0^{M\infty}$ faster. 

\section*{Appendix}

{\bf Proof of Proposition~\ref{prop-p-l-u}}. Note (\ref{p-ci-l-u}) and $p_B(x,n,p_0)=p_B(n-x,n,1-p_0)$. So,
\begin{eqnarray*} 
	h_p(n-x,1-p_0)
	&=& \sum_{\{n-y: T_p(n-y,p_0)\leq T_p(x,p_0)\}} p_B(n-y,n,p_0) \\
	&=& \sum_{\{y':T_p(y',p_0)\leq T_p(x,p_0)\}} p_B(y',n,p_0)=h_p(x,p_0). 
\end{eqnarray*}
Therefore, $
h_p(x,U_p(x))=h_p(n-x,1-U_p(x)),$
establishing (\ref{p-l-u}). 
~\raisebox{.5ex}{\fbox{}}
\medskip

{\bf Proof of Proposition~\ref{prop-diff-l-u}}. Note (\ref{diff-td}), $D(-d_0)=1-D(d_0)$ and $p_B(x,n,p)=p_B(n-x,n,1-p)$. We have
\begin{eqnarray*} 
	&&	h_d(n_1-x,n_2-y,-d_0)\\
	&=& \sup_{p_2\in D(-d_0)} \sum_{\{(u,v)\in S_d: T_d(u,v,-d_0)\leq T_d(n_1-x,n_2-y,-d_0)\}} p_B(u,n_1,p_2-d_0)\\ && \hspace{3in}*p_B(v,n_2,p_2) \\
	&=& \sup_{p_2\in D(-d_0)} \sum_{\{(u,v)\in S_d: T_d(u,v,-d_0)\leq T_d(x,y,d_0)\}} p_B(n_1-u,n_1,1-p_2+d_0)\\ && \hspace{3in}*p_B(n_2-v,n_2,1-p_2) \\
	&=& \sup_{p_2'\in D(d_0)} \sum_{\{(u',v')\in S_d: T_d(u',v',d_0)\leq T_d(x,y,d_0)\}} p_B(u',n_1,p_2'+d_0)p_B(v',n_2,p_2') \\
	&=&h(x,y,d_0). 
\end{eqnarray*}
Therefore, $h_d(x,y,U(x,y))=h_d(n_1-x,n_2-y,-U(x,y)),$	establishing (\ref{diff-ci-l-u}). 
~\raisebox{.5ex}{\fbox{}}
\medskip

{\bf Proof of Theorem~\ref{thm-ci-modi}.} Due to the definition of $h_2(\underline{x},\theta_0)$ in (\ref{h-2}), it is clear that  $h_2(\underline{x},\theta_0)$ is a valid p-value, see, for example, Casella and Berger (2002, p. 397).  Thus, $A_2(\theta_0)$ in (\ref{test-ci-2}) is the acceptance region of a level-$\alpha$ test. It implies that interval $C_0^M(\underline{x})$ is of level $1-\alpha$.
~\raisebox{.5ex}{\fbox{}}
\medskip	

{\bf Proof of Theorem~\ref{thm-ci-refine}.} We only need to prove the first claim because the second claim follows Theorem~\ref{thm-ci-modi} and the first claim. To prove the first claim, it suffices to show 
$h_2(\underline{x}, \theta_0)\leq \alpha$ for any $\theta_0 \not \in C_0(\underline{x})$.  Let $Cover_{C_0}(\theta,\underline{\eta})$ be the coverage probability function of interval $C_0(\underline{X})$ with an infimum at least $1-\alpha$ over the entire parameter space $H$.   

First consider the case of  $\theta_0< L_0(\underline{x})$. For the null hypothesis $H_0:\theta=\theta_0$,
\begin{eqnarray*}
	h_2(\underline{x}, \theta_0)&=&\sup_{H_0}P(T_2(\underline{X},\theta_0)\leq T_2(\underline{x},\theta_0) )\\ &=&\sup_{H_0}P(T_2(\underline{X},\theta_0)\leq \min\{\theta_0-L_0(\underline{x}), U_0(\underline{x})-\theta_0 \} )\\
	&\stackrel{\theta_0< L_0(\underline{x})}{\leq}
	&\sup_{H_0}P(T_2(\underline{X},\theta_0)<0 )\\
	&=&\sup_{H_0}(1-P(T_2(\underline{X},\theta_0)\geq 0 ))\\
	&\stackrel{(due\,\ to \,\ \ref{t2-c0})}{=}&\sup_{H_0}(1-P(\theta_0\in C_0(\underline{X})))\\
	&=&1- \inf_{H_0} Cover_{C_0}(\theta_0,\underline{\eta})\\
	&\leq& 1- \inf_H Cover_{C_0}(\theta,\underline{\eta})\leq \alpha.
\end{eqnarray*}
When $\theta_0>U_0(\underline{x})$, we establish $h_2(\underline{x},\theta_0)\leq \alpha$ similarly. ~\raisebox{.5ex}{\fbox{}}

\medskip
{\bf Proof of Theorem~\ref{thm-ci-refine-inf}.} We only prove part ii) as the other claims are straightforward. 
Let $Cover_C(\theta,\underline{\eta})$ be the coverage probability function for an interval $C(\underline{X})$. 
Note that the indicator functions satisfy 
$$I_{C_0^{M\infty}(\underline{x})}(\theta)=\lim_{k\rightarrow +\infty} I_{C_0^{Mk}(\underline{x})}(\theta), \,\ \forall \underline{x}$$
because $C_0^{Mk}(\underline{x})$ is nonincreasing and note that
$Cover_{C_0^{Mk}}(\theta,\underline{\eta})\geq 1-\alpha$ for any $(\theta,\underline{\eta})$ because each interval $C_0^{Mk}(\underline{X})$ is of level $1-\alpha$. 
Then, following the Dominated Convergence Theorem
\begin{eqnarray*}
	Cover_{C_0^{M\infty}}(\theta,\underline{\eta})&=& E_{(\theta,\underline{\eta})}[I_{C_0^{M\infty}(\underline{x})}(\theta)] = \lim_{k \rightarrow +\infty}E_{(\theta,\underline{\eta})}[I_{C_0^{Mk}(\underline{x})}(\theta)]\\ &=&\lim_{k \rightarrow +\infty}	Cover_{C_0^{Mk}}(\theta,\underline{\eta})
	\geq 1-\alpha,
\end{eqnarray*}
which completes the proof.
~\raisebox{.5ex}{\fbox{}}

{\bf Proof of Theorem~\ref{thm-ci-1l}}. Part i) is similar to the proof of Theorem~\ref{thm-ci-modi}. Part iii) is similar to the proof of Theorem 4 in Wang (2010) and is skipped. Part ii) follows part iii). 
~\raisebox{.5ex}{\fbox{}}
\medskip

\setlength{\parindent}{-1.5em}

\begin{Large}
	{\bf References}
\end{Large}
\bigskip

Agresti, A. (2013). {\it Categorical Data Analysis}, 3rd edition. John Wiley \& Sons.

Agresti, A. and Min, Y. (2001). ``On Small-sample Confidence Intervals for Parameters in Discrete Distributions,'' {\it Biometrics} 57, 963-971.
 
Arbuthnot, J. (1710). ``An Argument for Divine Providence, Taken from the Constant Regularity Observed in the Births of Both Sexes,'' {\it Philosophical Transactions of the Royal Society of London}, 27 (325-336). 186-190.


Bentur, L., Lapidot, M., Livnat, G., Hakim, F., Lidroneta-Katz, C., Porat,
I., Vilozni, D., and Elhasid, R. (2009) ``Airway Reactivity in Children Before and After Stem Cell Transplantation,'' {\it Pediatric Pulmonology}, 44, 845-850.

Blaker, H. (2000). ``Confidence Curves and Improved Exact Confidence Intervals for Discrete Distributions,'' {\it Canadian Journal of Statistics}, 28, 783-798.

Bonett, D. G. and Price, R. M. (2012). ``Adjusted Wald Confidence Interval for a Difference of Binomial Proportions Based on Paired Data,'' {\it Journal of Educational and	Behavioral Statistics}, 37(4), 479-488.

Brown, L. D., Cai, T. T. and DasGupta, A. (2001). ``Interval Estimation for a Binomial Proportion,'' {\it Statistical Science}, 16, 101-133.

Buehler, R. J. (1957). ``Confidence Intervals for the Product of Two Binomial Parameters,'' {\it Journal of the American Statistical Association}, 52, 482-493.
	
Casella, G. (1986), ``Refining Binomial Confidence Intervals,'' 
	\textit{Canadian Journal of Statistics}, 14, 113-129.
		
Casella, G., and Berger, R. L. (2002), \textit{Statistical Inference,} 2nd ed. Duxbury Press, Pacific Grove, CA.
	
Clopper, C. J., and Pearson, E. S. (1934). ``The Use of Confidence or Fiducial Limits in the Case of the Binomial,'' \textit{Biometrika}, 26, 404-413.

Essenberg, J. M. (1952). ``Cigarette Smoke and the Incidence of Primary Neoplasm of the Lung in Albino Mice,'' {\it Science}, 116, 561-562.

Fagerland, M. W., Lydersenb, S. and Laakec, P. (2013). ``Recommended Tests and Confidence Intervals for Paired Binomial Proportions,'' {\it Statistics in Medicine}, 33(16), 2850-2875.

Fay, M. P. (2010). ``Two-sided Exact Tests and Matching Confidence Intervals for Discrete Data,'' {\it The R Journal}, 2(1), 53-58.

Fisher, R. (1925). {\it Statistical Methods for Research Workers}. Edinburgh, Scotland: Oliver \& Boyd.

Huwang, L. (1995). ``A Note on the Accuracy of an Approximate Interval for the Binomial Parameter.'' \textit{Statistics and Probability	Letters}, 24, 177-180.

Pearson, K. (1900). ``On the Criterion That a Given System of Deviations from the Probable in the Case of a Correlated System of Variables is Such That it can be Reasonably Supposed to Have Arisen from Random Sampling,'' {\it Philosophical Magazine Series} 5, 50 (302), 157-175.

Newcombe, R. G. (1998). ``Improved Confidence Intervals for the Difference Between Binomial Proportions Based on Paired Data,'' {\it Statistics in Medicine} 17, 2635-2650.

Shan, G. and Wang, W. (2013). ``ExactCIdiff: An R Package for Computing Exact
	Confidence Intervals for the Difference of Two Proportions,'' {\it The R Journal}, 5, 62-70.
	
Spjotvoll, E. (1983). ``Preference Functions,'' {\it In A Festschrift for Erich L. Lehmann},	P. J. Bickel, K. Doksum, and J. L. Hodges Jr, (Eds.) 409-432.

Stigler, S. M. (1986). {\it The History of Statistics: The Measurement of Uncertainty Before 1900}. Harvard University Press.

Tango, T. (1998). ``Equivalence Test and Confidence Interval for the Difference in Proportions for the Paired-sample Design,'' {\it Statistics in Medicine}, 17, 891-908.

Wang, H. (2007). ``Exact Confidence Coefficients of Confidence Intervals for a Binomial Proportion,'' {\it Statistica Sinica}, 17, 361-368.

Wang, W. (2010). ``On Construction of the Smallest One-sided Confidence Interval for the Difference of Two Proportions,'' {\it The Annals of Statistics}, 38, 1227-1243.
		
Wang, W. (2014). ``An Iterative Construction of Confidence Interval for a Proportion,''  \textit{Statistica Sinica}, 24, 1389-1410.
		
Wang, W. and Zhang, Z. (2014). ``Asymptotic Infimum Coverage Probability for
Interval Estimation of Proportions,'' {\it Metrika}, 77, 635-646.

Wang, W. (2018). ``A `paradox' in Confidence Interval Construction Using Sufficient Statistics,'' {\it The American Statistician}, 72(4), 315-320.
		
Wilson, E. B. (1927). ``Probable Inference, the Law of Succession, and Statistical Inference,''  \textit{Journal of the American Statistical Association}, 22, 209-212.
		

\end{document}